\newif\ifconfver
\confverfalse      

\newif\ifcutshort      
\cutshorttrue

\newif\ifcutshortlvltwo  
\cutshortlvltwofalse

\ifconfver
    \documentclass[journal]{IEEEtran}
\else
    \documentclass[11pt,draftcls,onecolumn]{IEEEtran}
\fi

\usepackage{psfrag,cite,amsmath,graphicx,epsfig,latexsym,subfigure,color}

\usepackage[scriptsize,bf]{caption}
\usepackage{color}
\usepackage{amssymb}
\usepackage{setspace}
\usepackage{array}
\usepackage{booktabs}

\def\bn{\hfill \\ \smallskip\noindent}

\def\prox{\mbox{prox}}

\newcommand{\beq}{\begin{equation}}
\newcommand{\eeq}{\end{equation}}
\newcommand{\st}{{\rm s.t.}}
\newcommand{\ag}{{\rm ag}}
\newcommand{\md}{{\rm md}}

\newcommand{\cG}{{\mbox{$\mathcal{G}$}}}
\newcommand{\cA}{{\mbox{$\mathcal{A}$}}}
\newcommand{\cV}{{\mbox{$\mathcal{V}$}}}

\newcommand{\barw}{{\mbox{$\bar{w}$}}}

\newcommand{\tP}{{\mbox{$\widetilde{P}$}}}

\begin{document}
\def\pn {\par\smallskip\noindent}
\def \bn {\hfill \\ \smallskip\noindent}
\newcommand{\fs}{f_1,\ldots,f_s}
\newcommand{\f}{\vec{f}}
\newcommand{\hf}{\hat{f}}
\newcommand{\hx}{\hat{x}}
\newcommand{\hy}{\hat{y}}
\newcommand{\hz}{\hat{z}}
\newcommand{\hw}{\hat{w}}
\newcommand{\tw}{\tilde{w}}
\newcommand{\hlambda}{\hat{\lambda}}
\newcommand{\hbeta}{\hat{\beta}}
\newcommand{\tG}{\widetilde{G}}
\newcommand{\tg}{\widetilde{g}}
\newcommand{\barhx}{\bar{\hat{x}}}
\newcommand{\vecx}{x_1,\ldots,x_m}
\newcommand{\xoy}{x\rightarrow y}
\newcommand{\barx}{{\bar x}}
\newcommand{\bary}{{\bar y}}
\newcommand{\hrho}{\widehat{\rho}}
\newtheorem{theorem}{\noindent\bf Theorem}[section]
\newtheorem{lemma}{\noindent\bf Lemma}[section]
\newtheorem{corollary}{\noindent\bf Corollary}[section]
\newtheorem{proposition}{\noindent\bf Proposition}[section]
\newtheorem{definition}{\noindent\bf Definition}[section]
\newtheorem{claim}{\noindent\bf Claim}[section]
\newtheorem{remark}{\noindent\bf Remark}[section]

\def\br{\break}
\def\smskip{\par\vskip 5 pt}
\def\proof{\bn {\bf Proof.} }
\def\QED{\hfill{\bf Q.E.D.}\smskip}
\def\qed{\quad{\bf q.e.d.}\smskip}
\newcommand{\cE}{\mathcal{E}}
\newcommand{\cM}{\mathcal{M}}
\newcommand{\cN}{\mathcal{N}}
\newcommand{\cJ}{\mathcal{J}}
\newcommand{\cT}{\mathcal{T}}
\newcommand{\bx}{\mathbf{x}}
\newcommand{\bp}{\mathbf{p}}
\newcommand{\bz}{\mathbf{z}}
\newcommand{\cF}{\mathcal{F}}

\newcommand{\blue}{\color{blue}}
\newcommand{\red}{\color{red}}

\title{Stochastic Proximal Gradient Consensus Over Random Networks}
\author{{Mingyi Hong and Tsung-Hui Chang}
\thanks{M. Hong is with the Department of Industrial and Manufacturing Systems Engineering, Iowa State University. Email: \texttt{mingyi@iastate.edu}. T.-H. Chang is with the School of Science \& Engineering, The Chinese University of Hong Kong, Shenzhen, E-mail: \texttt{tsunghui.chang@ieee.org}}\thanks{The conference version of this work appears in ICASSP 2016 \cite{hong_chang_15_icassp}.}
}
\vspace{-2cm}
\maketitle
\vspace{-2cm}
\begin{abstract}
We consider solving a convex optimization problem with possibly stochastic gradient, and over a randomly time-varying multi-agent network.   Each agent has access to some  local objective function, and it only has unbiased estimates of the gradients of the smooth component.
We develop a dynamic stochastic proximal-gradient consensus (DySPGC) algorithm,  with the following key features: {\it i)} it works for both the static and certain randomly time-varying networks;
{\it ii)} it allows the agents to utilize either the exact or stochastic gradient information; {\it iii)} it is convergent with provable rate. In particular, the proposed algorithm converges to a global optimal solution, with a rate of $\mathcal{O}(1/r)$ [resp. $\mathcal{O}(1/\sqrt{r})$] when the exact (resp. stochastic) gradient is available, where $r$ is the iteration counter.
Interestingly, the developed algorithm {establishes a close connection among} a number of (seemingly unrelated) distributed  algorithms, such as the EXTRA (Shi {\it et al.} 2014), the PG-EXTRA (Shi {\it et al.} 2015), the IC/IDC-ADMM (Chang {\it et al.} 2014),  the DLM (Ling {\it et al.} 2015) and the classical distributed subgradient method.  

\end{abstract}
\vspace{-0.3cm}
\section{Introduction}
\subsection{The Global Consensus Problem}

Consider the following classical problem
\begin{align}\label{eq:global:consensus}
\min_{y} \; f(y):=\sum_{i=1}^{N} f_i(y),
\end{align}
where {$f_i:\mathbb{R}^M\to \mathbb{R}$} is a convex and possibly nonsmooth function, for $i=1,2,\cdots, N$. Consider a collection of $N$ agents connected by a network defined by an {\it undirected} graph $\cG=\{\mathcal{V}, \mathcal{E}\}$, with $|\mathcal{V}|=N$ vertices and $|\mathcal{E}|=E$ edges. Each agent $i\in \mathcal{V}$ can communicate with its immediate neighbors, and it {only has information about its local function $f_i$, but not the function $f_j$ of any other agent $j\ne i$}. 

This problem has found applications in various domains such as distributed consensus \cite{tsitsiklis84thesis,Xiao:2007:DAC:1222667.1222952}, distributed and parallel machine learning \cite{Forero11,mateos10, shalev14proximaldual} and distributed signal processing \cite{Schizas08,zhu10}; see \cite{Giannakis15} for a recent survey. The key research question is: How to  compute an optimal solution of \eqref{eq:global:consensus}, through a distributed process where each agent only utilizes local gradient information about the objective.


Let each agent $i$ keep a local copy of $y$, say $y_i$. The well-known distributed subgradient (DSG) method \cite{Nedic09subgradient} is given by
\begin{align}
{y_i^{r+1} = {\sum_{j=1}^{N}}w^r_{ij} y^r_j - \gamma^r d^r_i, \; \forall~i\in\cV,\label{eq:distributed:subgradient:method}}
\end{align}
where $r$ denotes the iteration counter; $d^r_i\in\partial f_i(y^r_i)$ denotes a subgradient of the local function $f_i$ evaluated at $y^r_i$; $w^r_{ij}\ge 0$ denotes the weight for the link {$e_{ij}\in \mathcal{E}$} at iteration $r$; and $\gamma^r> 0$ denotes some stepsize parameter. Let $\bar{y}^r_i:=\frac{1}{r}\sum_{t=1}^{r}y_i^t$.

The convergence of the DSG iteration \eqref{eq:distributed:subgradient:method} was first analyzed in \cite{Nedic09subgradient} by Nedi\'{c} and Ozdaglar. It was shown that if the { subdifferential} is bounded, and that the weights $\{w_{ij}\}$ and the graph $\cG$ satisfy certain regularity assumptions, then each $\bar{y}^r_i$ converges to a neighborhood of the optimal solution (resp. the exact optimal solution) if $\gamma^r$ is a constant (resp. a diminishing sequence). As a special case, when $f(x)\equiv 0$ (only the consensus among the agents is sought for), then the convergence of the iteration \eqref{eq:distributed:subgradient:method} was first studied by Tsitsiklis \cite{tsitsiklis84thesis}. 
The DSG iteration has been extended to scenarios where there is a local constraint for each agent \cite{Nedic10constrained}, or the messages exchanged among the agents are quantized \cite{Nedic08sub_quantize}, or the communication among the agents is noisy \cite{Srivastava11}. Also see \cite{shi14extra,shi15pgextra,chen12thesis,Jakovetic14,Duchi12, Nedic15Time-Varying} for other related methods for solving \eqref{eq:global:consensus}. 

The rate of convergence analysis of the DSG-type method has been a central research issue. In its most general form, it is known that when appropriate diminishing stepsizes are chosen, DSG converges with a rate of $\mathcal{O}(\ln(r)/\sqrt{r})$ {in terms of the differences between the local objective functions and the optimal objective function} \cite{chen12thesis}, for both static and time-varying networks.  Duchi {\it et al.} propose a distributed dual-averaging algorithm and show that it converges with a rate of $\mathcal{O}(\ln(r)/\sqrt{r})$.  Jakovetic {\it et al.} \cite{Jakovetic14} show that when the objective has Lipschitz continuous and bounded gradient, and when the graph is static, it is possible to accelerate the DSG to achieve an $\mathcal{O}(1/r^2)$ rate, but at the expense of solving more complicated subproblems, each of which involves multiple rounds of communication and computation. If only simple computation/communication steps are performed, the rate becomes $\mathcal{O}(\ln(r)/r)$. A related acceleration scheme has also been proposed in \cite{chen12thesis}, which further works for time-varying $B$-connected graphs \footnote{A $B$-connected graph is a time-varying graph in which at each iteration the graph is not necessarily connected, but the union of the graphs across every $B>0$ consecutive iterations is connected.}. Under the smoothness assumption on $f$, Shi {\it et al.} \cite{shi14extra} propose an interesting algorithm called EXTRA, which adds certain {\it error-correction} terms to the DSG \eqref{eq:distributed:subgradient:method}. By adding such {correction}, EXTRA uses constant stepsize and achieves an $\mathcal{O}(1/r)$ rate for smooth convex problem and linear convergence for certain smooth strongly convex problems. This method has also been generalized to solve nonsmooth problems \cite{shi15pgextra}, but both algorithms in \cite{shi15pgextra,shi14extra} can only work for static networks. {Other recent developments can be found in \cite{Mokhtari16DSA,Mokhtari16linear} and the references therein.}

Another popular approach for distributed optimization is to use the alternating direction method of multipliers (ADMM) \cite{Bertsekas_Book_Distr,BoydADMMsurvey2011,Glow84}. Applying the ADMM to distributed optimization has been first suggested in \cite{Bertsekas_Book_Distr}, and subsequently popularized in \cite{Schizas08,BoydADMMsurvey2011}. The $\mathcal{O}(1/r)$ sublinear rate of convergence for decentralized consensus ADMM (C-ADMM) has been shown by Wei and Ozdaglar \cite{Wei13}, where it is assumed that the underlying graph is generated according to certain stochastic mechanism. When the problem is smooth, the linear convergence of C-ADMM is shown in \cite{shi14linear}. Recently a broadcast based C-ADMM has been proposed in \cite{Makhdoumi14}. However the C-ADMM usually requires solving local optimization problems exactly (cf. \cite{Mota13,Schizas08,zhu10,BoydADMMsurvey2011,Erseghe11fast_consensus,Makhdoumi14}), which can be expensive in certain applications. This requirement has been relaxed by two recent works  \cite{chang14distributed} and \cite{Ling15decentralize}. In particular, Chang {\it et al.} \cite{chang14distributed} develop an inexact C-ADMM (IC-ADMM) algorithm 
which uses a simple (proximal) gradient step at each ADMM iteration. Ling {\it et al.} \cite{Ling15decentralize} also propose to replace the exact minimization by certain proximal gradient steps. {While we are finalizing the paper, we were made aware of an independent work \cite{aybat16} that also proposes linearlized ADMM method for consensus composite optimization. In particular, for a static network, convergence rates of $\mathcal{O}(1/{r})$ and $\mathcal{O}(1/\sqrt{r})$ are shown for certain deterministic and stochastic linearlized distributed ADMM.} Recently, Hong {\it et al.} \cite{hong14nonconvex_admm} show that the ADMM-based method (with exact or inexact update) can  be used to solve certain {\it nonconvex} global consensus problem, with a convergence rate of $\mathcal{O}(1/\sqrt{r})$. 

There has been a few works that design distributed optimization algorithms in the primal-dual perspective. For example, \cite{Pesquet15, Bianchi14, combettes15} propose random coordinate primal-dual algorithms, with possible applications in distributed and  asynchronous optimization.  However no convergence rate has been provided. In \cite{Jakovetic14b}, the authors propose an augmented Lagrangian method based algorithm for distributed optimization and analyzed its linear convergence, but the algorithm and analysis only works for smooth and strongly convex problems. Further, the algorithm has double loops, and requires some global knowledge about the objective function and the underlying graph. These requirements can be restrictive in practical applications.  {Reference \cite{Mokhtari16DSA} develops a primal-dual algorithm in which each agent is updated by performing local stochastic averaging gradients. }

Below we provide a high level comparison of the DSG-based and ADMM-based algorithms.
\begin{itemize}
	\item \textsf{(Problem types)} The DSG can solve convex problems with only subgradient information about the objective, while to our best knowledge the ADMM does not directly work for this case. 
	\item \textsf{(Gradient Information)} The DSG only needs (stochastic) subgradients of the objective \cite{Ram2010_stoc}, {while the ADMM usually requires subproblems to have some nice structures so that they can be solved in closed-form \cite{Schizas08,BoydADMMsurvey2011}}. 
	\item \textsf{(Convergence rates)} When the objective function $f$ has certain additional structures (e.g., smooth or a smooth plus a simple nonsmooth function), {the distributed ADMM generally converges faster in practice} (which has a convergence rate of $\mathcal{O}(1/r)$) than its DSG counterpart (which has a convergence rate of $\mathcal{O}(\ln(r)/\sqrt{r})$). {Theoretically, it is possible to modify the iteration of the DSG algorithm to improve its rate to $\mathcal{O}(1/r)$; see recent developments in  \cite{shi14extra,Jakovetic14,chen12thesis}}.
	\item \textsf{(Network structures)} The DSG generally works when the underlying network is time-varying and follows the so-called $B$-connected structure \cite{Nedic09}. However the ADMM-based method only works for static network, except for the recent variants proposed in \cite{Wei13,chang14}, both of which work for certain randomized networks.
\end{itemize}

%
%
%
%
%

\vspace{-0.3cm}
\subsection{Contribution of This Work}

In this work, we consider the following structured version of the global consensus problem \eqref{eq:global:consensus}{
\begin{align}\label{eq:consensus}
\min_{y}\; f(y):=\sum_{i=1}^{N}f_i(y)=\sum_{i=1}^{N} \left(g_i(y) +h_i(y)\right),
\end{align}}
\!\!where each $g_i: \mathbb{R}^M\to \mathbb{R}$ is a smooth convex function; each {$h_i: \mbox{dom}(h_i)\to \mathbb{R}\cup \{\infty\}$} is a convex possibly lower semi-continuous function, which covers the indicator function for closed convex sets as  a special case. 

We propose an ADMM based method, named dynamic stochastic proximal-gradient consensus (DySPGC), that has the following key features:

\noindent \textbullet~When only an unbiased estimate of each $\nabla g_i$ is known, the algorithm converges with a rate $\mathcal{O}(1/\sqrt{r})$;

\noindent \textbullet~When the exact $\nabla g_i$ is known, the rate becomes $\mathcal{O}(1/r)$;

\noindent \textbullet~The algorithm works for both the static and certain random time-varying networks.

What is more interesting is our insight on the connection between the C-ADMM-type methods and a few DSG-type methods. In particular, we show that the EXTRA/PG-EXTRA \cite{shi14extra,shi15pgextra}, despite being posed as error-corrected DSGs, can be viewed as special cases of the proposed DySPGC (for static network with symmetric weights and exact gradients). This observation explains the relative fast practical convergence performance of these two algorithms compared with the  DSG [for structured problems \eqref{eq:consensus}]. Further, we also establish a close connection between the DSG \eqref{eq:distributed:subgradient:method} and the proposed DySPGC. Additionally our method generalizes other distributed ADMM-type methods such as the DLM \cite{Ling15decentralize} and the IC-ADMM \cite{chang14distributed}.


\vspace{-0.2cm}
\section{System Model}

\subsection{Problem Setup}
We  make the following blanket assumptions for  \eqref{eq:consensus}. 
\pn{\bf Assumption 1.}
\begin{enumerate}
\item The optimal solution set  of \eqref{eq:consensus}, denoted as $X^*\subseteq \mathbb{R}^{M}$ is nonempty; {The Slater condition holds for problem \eqref{eq:consensus};}
 \item Let $g(y):=\sum_{i=1}^{N} g_i(y)$, and  $h(y):=\sum_{i=1}^{N} h_i(y)$. The $h_i$'s {\it prox operators}, defined below
\begin{align}
\prox^{\beta}_{h_i}(u): = \min_{y}\; h_i(y) + \frac{\beta}{2}\|y-u\|^2,\label{eq:prox}
\end{align}
are easy to compute;
 \item Each $\nabla g_i$ is Lipschitz continuous (with constant $P_i>0$)
\begin{align}\label{eq:lip}
\begin{split}
\hspace{-0.3cm}\|\nabla g_i(y)-\nabla g_i(v)\|&\le P_i \|y-v\|, \; \forall~y, v \in \mbox{dom}(h).
\end{split}
\end{align}
\end{enumerate}
As have been mentioned in the introduction, we consider a collection of $N$ agents defined over a connected {\it undirected} graph $\cG=\{\mathcal{V}, \mathcal{E}\}$, with $|\mathcal{V}|=N$ vertices and $|\mathcal{E}|=E$ edges. Define a companion {\it symmetric} directed graph given by $\cG_d=\{\mathcal{V}, \mathcal{A}, W\}$, where $\mathcal{A}$ is a set of directed arcs with $|\mathcal{A}|=2E$, and for every edge in $\mathcal{E}$ which connects nodes $i,j$, we have  both $e_{ij}, e_{ji}\in\mathcal{A}$. {Note that using a companion graph to represent the original graph $\cG$ is conventional in the consensus ADMM literature; see e.g., \cite{chang14distributed} and \cite{Ling15decentralize}. It helps to simplify the definition of the consensus constraint (to be provided shortly).} 

Let us use $\cN_i$ to denote the neighborhood of node $i$, i.e.,
\begin{align}
\cN_i := \{ j\mid e_{ij}\in \cA\}. 
\label{eq:neighbor}
\end{align}
{In distributed optimization, we also need a  {\it weight matrix} 
$W\in\mathbb{R}_{+}^{N\times N}$  whose coefficients are used by the agents to combine their neighboring messages (see e.g., \eqref{eq:distributed:subgradient:method}).}
Generally we will assume that the weight matrix $W$ satisfies the following two conditions:
\begin{enumerate}
\item $W$ is a row stochastic matrix, i.e., $\{W[i,j]\ge0\}$, $\sum_{j}W[i,j]=1, \; \forall~i$;
\item The diagonal elements of $W$ are all positive, and its off-diagonal elements all satisfy
\begin{align}\label{eq:w:positive}
W[i,j]>0, \; \mbox{if} \; e_{ij}\in\cA, \quad W[i,j]=0, \; \mbox{otherwise}.
\end{align}
\end{enumerate}
Later we will provide explicit expressions for $W$.


Consider an equivalent reformulation of problem \eqref{eq:consensus} (equivalent when $\mathcal{G}$ is connected) 
\begin{align}\label{eq:consensus:admm}
\begin{split}
\min_{\{x_i\},\{z_{ij}\}}&\quad f(x):=\sum_{i=1}^{N}\left( g_i(x_i) +h_i(x_i)\right),\\
\st &\quad x_i = z_{ij}, \quad x_j = z_{ij}, \; \forall~(i,j)\in\mathcal{A},
\end{split}
\end{align}
where we have introduced $N$ auxiliary variables $\{x_i\in\mathbb{R}^M\}$, and $2E$ auxiliary variables $\{z_{ij}\in\mathbb{R}^M\}$. Define $x: = \{x_i\}\in\mathbb{R}^{NM\times 1}$, and $z:=\{z_{ij}\}\in\mathbb{R}^{2EM}$.
To compactly represent the constraint set of problem \eqref{eq:consensus:admm}, let us define the following two matrices
\begin{align}\label{eq:A}
A := \left[\begin{array}{l}
A_1\\
A_2
\end{array}\right], \; B := \left[\begin{array}{l}
-I_{2EM}\\
-I_{2EM}
\end{array}\right],
\end{align}
where $A_1, A_2\in\mathbb{R}^{2EM\times NM}$, and each of them is composed of $2E\times N$ blocks of $M\times M$ matrices. If $e_{ij}\in\mathcal{A}$ and $z_{ij}$ is the $q$th block of $z$, then the $(q,i)$th block of $A_1$ and the $(q,j)$th block of $A_2$ are both $I_M$, an $M\times M$ identity matrix; otherwise, the corresponding block is an $M\times M$ zero matrix $0_{M}$. Note that the matrix $B$ stacks two identity matrices because each link variable $z_{ij}$ only appears once in the constraint. 
Using the above matrix notation, problem \eqref{eq:consensus:admm} is equivalent to the following problem \cite{Schizas08,Ling15decentralize,Giannakis15,chang14distributed}
\begin{align}\label{eq:consensus:admm:matrix}\tag{P}
\min_{x,z}&\quad f(x):=\sum_{i=1}^{N} \left(g_i(x_i) +h_i(x_i)\right)\\
\st & \quad  Ax+Bz = 0. \notag
\end{align}

\vspace{-0.3cm}
\subsection{Randomly Time-Varying Graph Structure}\label{sub:random}
We assume that the edges of the graph $\cG$ are activated according to certain randomly time-varying patterns. 
To describe such random pattern, at a given time $r$, define a new graph $\cG^{r}=\{\cV^{r}, \cE^{r}\}$, and its companion graph $\cG_d^{r} = \{\cV^{r}, \cA^{r}, W^{r}\}$ where $\cV^{r}\subseteq \cV$, $\cE^{r}\subseteq\cE$ and $\cA^{r}\subseteq \cA$, and each weight matrix $W^{r}$ is a stochastic matrix satisfying \eqref{eq:w:positive}. Again $\cG_d^{r}$ is symmetric, meaning if $e$ connects nodes $i$ and $j$ with $e\in\cE^{r}$, then $e_{ij}, e_{ji}\in\cA^{r}$. 
The precise specification of the random graphs $\{\cG_d^r\}$ and $\{\cG^r\}$  is given below \cite{Wei13,Srivastava11,chang14,Duchi12}.
\begin{definition}(Randomly Activated Graph)\label{def:graph}
At each time $r$, each link pair $(i,j), (j,i)\in \cA$ has a probability $p_{ij}=p_{ji}\in(0,1]$ of being active. The set of active nodes $\cV^{r}$ is given by:
    $$\cV^{r}=\{i\mid \exists e_{ij}\in \cA^{r}, \; \forall~j\in\cV\}.$$
Effectively at each time $r$  a node $i\in \cV$ has a probability $\alpha_i>0$ of being active, while such $\alpha_i$ is a function of $\{p_{ij}\mid j\in\cN_i\}$. Let us collect these probabilities and define
\begin{align}\label{eq:Psi}
\Psi=\mbox{diag}\{\alpha_i\}\in\mathbb{R}^{N\times N}, \quad \Phi=\mbox{diag}\{p_{ij}\}\in\mathbb{R}^{2E\times 2E}
\end{align}
{where $\mbox{diag}\{\alpha_i\}$ represents a diagonal matrix whose diagonal entries are the elements in the set $\{\alpha_i\}$.
Further, assume that $\cG_d$ is strongly connected, and realizations of the graphs $\cG_d^{r}$ and $\cG_d^t$ are independent and identically distributed across all $r\ne t$}. \hfill $\blacksquare$
    \end{definition}

In practice, the random network pattern can be used to model communication and/or node failures \cite{Wei13,Srivastava11,chang14,Duchi12}. It is the stochastic variant of the so-called $B$-strongly connected network which has been widely considered in the literature, under very different context \cite{zhu10,tsitsiklis84thesis, Nedic09,Nedic15Time-Varying}.
The connection between such randomly generated graph and popular communication protocols such as the gossip protocol and asynchronous protocols has been explored in \cite{Srivastava11,Duchi12, Wei13}. Note the graph $\cG$ is required to be connected, but $\cG^r$'s are not necessarily so. At a given iteration $r$, we can define the neighborhood $\cN_i^{r}$ for each node $i$ similarly as in \eqref{eq:neighbor}, and define the matrices $A^{r}$ and $B^{r}$ similarly as in \eqref{eq:A}, making all quantities conforming to the instantaneous graph structure. 

\vspace{-0.3cm}
\subsection{The Gradient Information}
Define the gradient of the smooth part of the objective as $G(x):=[\nabla g_1(x_1);\cdots; \nabla g_N(x_N)]$.
In this work, we will consider situations in which only an estimate of $\nabla g_i(x_i)$, denoted by $\tg_i( x_i, \xi_i)$, is available for each agent $i$. \
In this case, the estimate $\tg_i( x_i, \xi_i)$ will satisfy the following
\begin{align}
&\mathbb{E}[\tg_i(x_i,\xi_i)] = \nabla g_i(x_i),\; \mathbb{E}\left[\|\tg_i(x_i,\xi_i)-\nabla g_i(x_i)\|^2\right] \le \sigma^2,\label{eq:gradient}
\end{align}
{where each $\xi_i$ is a random variable following an unknown distribution, and $\xi_i$, $\xi_j$ are not necessarily independent for any $i\ne j$. Further, when time is involved (cf. Section \ref{sub:random}), we will assume $\xi_i$ to be independent over time.  Each $\tg_i( x_i, \xi_i)$ is assumed to be a measurable function; the constant $\sigma^2$ represents the maximum expected deviation of the gradient estimate.} 

\vspace{-0.3cm}
\section{The Proposed Algorithms}\label{sec:Algorithms}
Our proposed algorithms are based on the ADMM. {To describe the algorithm in its general form, let us first define a vector of positive {\it penalty constants} $\rho:=\{\rho_{ij}>0\mid e_{ij}\in\cA\}$, i.e., each $\rho_{ij}$ corresponds to a link variable $z_{ij}$. For a given graph $\cG_{d}$, we can construct a {\it diagonal} matrix $\Gamma\succeq 0$ by}
{\begin{align}\label{eq:gamma}
	\Gamma =\mbox{blkdg}[\Xi\otimes I_{M}, \Xi\otimes I_{M}]\in \mathbb{R}^{4EM\times 4EM},
	\end{align}
	where the notation {\it $\mbox{blkdg}$} represents taking the {\it block diagonal} operation;
	$\Xi\in\mathbb{R}^{2E\times 2E}$ is a diagonal matrix with $\Xi[q,q] = \rho_{ij}$ if link $(i,j)\in \cA$ and $z_{ij}$ is the $q$th block of $z$.}

Using the above definition, let us  write the augmented Lagrangian of \eqref{eq:consensus:admm:matrix}:
{\begin{align}\label{eq:augmented}
{ L}_{\Gamma}(x, z, \lambda) =  \sum_{i=1}^{N}f_i(x_i)+\langle \lambda, Ax+Bz \rangle+\frac{1}{2}\|Ax+Bz\|_{\Gamma}^2,
\end{align}}
\!\!where $\lambda\in\mathbb{R}^{4EM}$ is the dual variable corresponding to the inequality constraint $Ax+Bz=0$. Our definition of the augmented Lagrangian is slightly different from the standard definition due to the use of the graph-related positive definite penalty matrix $\Gamma\in\mathbb{S}_{++}^{4EM\times 4EM}$. Such modification turns out to be crucial in modeling some graph specific properties.

To proceed, we need the following definitions. For each $i\in\cV$ and some $\omega_i\ge 0$, define 
\begin{align}\label{eq:omega}
\Omega_i:= \omega_i I_M, \; \mbox{and}\; \Omega:= \mbox{blkdiag}\{\Omega_1,\cdots, \Omega_N\}, 
\end{align}
{where the latter matrix is a block diagonal matrix with diagonal blocks being $\Omega_1, \cdots, \Omega_N$}. Define the following matrices
\begin{align}
&M_{+}:= A^T_1+ A^T_2, \quad M_{-}:= A^T_1 - A^T_2. \label{eq:def:M}
\end{align}
{It can be verified that $\frac{1}{2}M_{-} M^T_{-}$ and $\frac{1}{2}M_{+} M^T_{+}$  represent the signed and signless graph Laplacian matrices, respectively (see, e.g.,  \cite[Section II]{Ling15decentralize} for detailed discussion on these matrices). }

{To illustrate various quantities related to the graphs, let us consider a simple graph with $3$ nodes and two edges connecting nodes $\{1, 2\}$ and nodes $\{2,3\}$. Suppose that $M=1$. In this case, $\mathcal{A}=\{(1,2), (2,3), (2,1), (3,2)\}$. Let us order the links as $(1,2), (2,3), (2,1), (3,2)$, then the matrices $A_1$ and $A_2$ are given below
	\begin{align*}
	A_1=\left[\begin{array}{lll}
	1 & 0& 0\\
	0& 1 & 0\\
	0& 1 & 0\\
	0 & 0& 1
	\end{array}\right], \quad A_2=\left[\begin{array}{lll}
	0 & 1& 0\\
	0& 0 & 1\\
	1& 0 & 0\\
	0 & 1& 0
	\end{array}\right].
	\end{align*}	
	The matrix $\Xi$ is given by
	\begin{align*}
	\Xi=\left[\begin{array}{llll}
	\rho_{12}& 0& 0& 0\\
	0& \rho_{23} & 0 &0\\
	0& 0 & \rho_{21}&0\\
	0 & 0& 0&\rho_{32}
	\end{array}\right].
	\end{align*}	
The matrices $M_{+}$ and $M_{-}$ are given by 	
	\begin{align*}
	M_{+}=\left[\begin{array}{llll}
	1& 0& 1& 0\\
	1& 1 & 1 &1\\
	0& 1 & 0&1
	\end{array}\right], \; 	M_{-}=\left[\begin{array}{llll}
	1& 0& -1& 0\\
	-1& 1 & 1 &-1\\
	0& -1 & 0&1
	\end{array}\right].	
	\end{align*}

}


\vspace{-0.3cm}
\subsection{The Proposed Algorithms}\label{sub:randomactivate:algorithm}
In this section we propose consensus algorithms over the randomly activated graphs.
To model the time-varying node activation pattern, let us define $\xi^{r+1}:=[\xi^{r+1}_1; \cdots; \xi^{r+1}_N]$, and define
$\tG^{r+1}(x^r,\xi^{r+1})\in\mathbb{R}^{MN}$ as a vector consisting of the gradients of the active component functions at time $r+1$, i.e., 
\begin{align*}
&\tG^{r+1}(x^r,\xi^{r+1}):=[a_1; \;a_2;\;\cdots\;;a_N]\\
&\mbox{with}\quad a_i = \left\{\begin{array}{ll}
\tg_i(x_i^{r},{\xi_i^{r+1}})\;\; &\mbox{if} \;\; i\in\cV^{r+1}\\
0 &\mbox{otherwise}
\end{array}\right..
\end{align*}
Define $h^{r+1}(x):=\sum_{i\in\cV^{r+1}}h_i(x_i)$. 
Let $\{\eta^{r}\ge 0\}$ denote a sequence of iteration-dependent parameters, whose values will be given shortly. 

Using these definitions, we present in the table below the proposed algorithm in its general form, named the {\it dynamic stochastic proximal-gradient consensus} (DySPGC) algorithm.

\hspace{-0.4cm}
\begin{center}
\vspace{-0.3cm}
\fbox{
\begin{minipage}{3.2in}
\smallskip
\centerline{\bf {Algorithm 1.  DySPGC Over Random  Graphs}}
\smallskip
{At iteration $0$, select $\lambda^0, z^0, x^0$ such that
$$B^T \lambda^0 = 0, \; z^0 = \frac{1}{2} M^T_{+} x^0.$$}

At each iteration $r+1$, update the variable blocks by:
{\small\begin{subequations}
\begin{align}
x^{r+1}& =\arg\min_x\;  \bigg\{ \left\langle \tG^{r+1}(x^{r},\xi^{r+1}), x-x^r\right\rangle + h^{r+1}(x) \notag \\
&~~~~~~~~~~~~~~+ \frac{1}{2}\left\|A^{r+1}x+B^{r+1}z^r+\Gamma^{-1}{\lambda^r}\right\|_{\Gamma}^2  \notag \\
&~~~~~~~~~~~~~~+\frac{1}{2}\|x-x^r\|^2_{\Omega+{\eta^{r+1}}I_{MN}} \bigg\} \label{eq:x:update3}\\
x^{r+1}_{i} & = x^r_i, \quad \mbox{if}~ i\notin\cV^{r+1}\label{eq:x:update3:unchange}\\
z^{r+1} & =\arg\min_z\; \frac{1}{2}\left\|A^{r+1} x^{r+1} + B^{r+1}z+\Gamma^{-1}{\lambda^r}\right\|_{\Gamma}^2\label{eq:z:update3}\\
z^{r+1}_{ij}&=z_{ij}^r,\quad\mbox{if}~e_{ij}\notin\cA^{r+1} \label{eq:z:update3:unchange}\\
\lambda^{r+1}& = \lambda^r +\Gamma\left( A^{r+1} x^{r+1}+ B^{r+1}z^{r+1}\right)\label{eq:lambda:update3}
\end{align}
\vspace{-0.5cm}
\end{subequations}}
\end{minipage}
}
\end{center}

Let us make a few comments about DySPGC. First, the penalty parameter used for the $x$-update for the proximal term $\|x-x^r\|^2$ is given by {$\Omega+ \eta^{r+1} I_{MN}$}. {Here {$\Omega$} is a fixed constant matrix defined in \eqref{eq:omega}; the iteration-dependent parameter $\eta^{r+1}$, when being chosen as an appropriate increasing sequence (to be specified in Theorem \ref{thm:time:variant:inexact:rate}),  is used to deal with the stochasticity in the gradient.} Second, when the gradients are precisely known, we can set $\eta^{r+1}=0$ for all $r$, in which case the $x$-update rule \eqref{eq:x:update3} becomes
\begin{align*}
x^{r+1}& =\arg\min_{x}\; \left\langle G^{r+1}(x^r), x-x^r\right\rangle + h^{r+1}(x) \nonumber\\\
&\quad+ \frac{1}{2}\left\|A^{r+1}x+B^{r+1}z^r+\Gamma^{-1}{\lambda^r}\right\|_{\Gamma}^2+\frac{1}{2}\|x-x^r\|^2_{\Omega}
\end{align*}
where $G^{r+1}(x^r)$ is defined similarly as $\tG^{r+1}(x^{r},\xi^{r+1})$ (with inexact gradients replaced by the exact gradients).

When we assume that the graph is static and the exact gradients are known, i.e., $\cG_d^{r}=\cG_d$ and $\tG^{r+1}(x^{r},\xi^{r+1})=G(x^r)$ for all $r$, then the DySPGC reduces to a simplified version named the  {\it proximal gradient consensus} (PGC) algorithm (see Algorithm 2).

{Let us compare Algorithms 1 and 2 with some  existing methods and pinpoint the main differences. First, PGC is a {\it proximal} version of the conventional C-ADMM \cite{Schizas08,BoydADMMsurvey2011,mateos10}, where we have used the second order approximation $\langle G(x^r), x-x^r\rangle+ \frac{1}{2}\|x-x^r\|^2_{\Omega}$ of the smooth function $g(x)$ in \eqref{eq:x:update:static} in the $x$-update, rather than exactly minimizing the augmented Lagrangian (as has been done in \cite{Schizas08,BoydADMMsurvey2011,mateos10}). Moreover, a matrix penalty $\Gamma$ is used instead of a scalar one; the latter has been popular in the existing ADMM-based methods. 
Later we will show that by using a matrix penalty, the parameters $\rho_{ij}$'s can be chosen by only using local information, making the algorithm better suited to distributed implementation. 
Second, the DySPGC is a {\it stochastic} version of the algorithms proposed in \cite{Wei13,chang14}, where we have used an iteration-dependent stochastic second order approximation
	$$\left\langle \tG^{r+1}(x^{r},\xi^{r+1}), x-x^r\right\rangle +\frac{1}{2}\|x-x^r\|^2_{\Omega+{\eta^{r+1}}I_{MN}}.$$ 
in the $x$-update step, rather than exactly minimizing the augmented Lagrangian

Detailed comparison with existing  algorithms will be provided in Section \ref{sec:compare}.}

\begin{center}
	\fbox{
		\begin{minipage}{3in}
			\smallskip
			\centerline{\bf Algorithm 2. PGC Over Static Graphs}
			\smallskip
			{At iteration $0$, select $\lambda^0, z^0, x^0$ such that
				$$B^T \lambda^0 = 0, \; z^0 = \frac{1}{2} M^T_{+} x^0.$$}
			
			At each iteration $r+1$, update the variable blocks by:
			\begin{subequations}
				\begin{align}
				x^{r+1}& =\arg\min_{x}\; \!\!\bigg\{\!\!\! \left\langle G(x^r), x-x^r\right\rangle \!+\! h(x)\!+\!\langle  \lambda^r, Ax+Bz^r\rangle\nonumber\\
				&~~~~~\quad + \frac{1}{2}\left\|Ax+Bz^r\right\|_{\Gamma}^2+\frac{1}{2}\|x-x^r\|^2_{\Omega} \bigg\} \label{eq:x:update:static}\\
				z^{r+1} & =\arg\min_{z}\; \frac{1}{2}\left\|A x^{r+1} +
				Bz+\Gamma^{-1}{\lambda^r}\right\|_{\Gamma}^2\label{eq:z:update:static}\\
				\lambda^{r+1}& = \lambda^r +\Gamma \left( A x^{r+1}+ Bz^{r+1}\right)\label{eq:lambda:update:static}
				\end{align}
			\end{subequations}
			
		\end{minipage}
	}
\end{center}

\subsection{Distributed Implementation}
Both algorithms proposed in the previous section can be implemented in a distributed manner, in which the information needed for updating each variable can be obtained from its immediate neighbors.  To see this, note that in the original formulation \eqref{eq:consensus:admm} each node $i$ is only coupled with its neighboring links $\{e_{ij}, e_{ji}\}_{j\in\mathcal{N}_i}$, and each link pair $e_{ij}, e_{ji}\in\cA$ is only related to its two neighboring nodes $\{i, j\} \in \cV$. Below we illustrate the distributed implementation of the PGC algorithm, as it takes a simple form.

To write the algorithm compactly, define the {\it stepsize parameter} $\beta_i$ as [with  $\hrho_{ij}:=1/2(\rho_{ij}+\rho_{ji})$]
\begin{align}\label{eq:beta}
\beta_i:=2 \bigg(\sum_{j\in\cN_i}\widehat{\rho}_{ij}+{ \omega_i}/{2}\bigg), \; \forall~i\in\cV.
\end{align}
Let us specialize the weight matrix $W\in\mathbb{R}^{N\times N}$ and define a new {\it stepsize matrix} ${\Upsilon}\in\mathbb{R}^{MN\times MN}$ as follows
\begin{align}\label{eq:W}
\hspace{-0.3cm} W[i,j]\hspace{-0.1cm}&=\hspace{-0.1cm}\left\{ \!\!\!\! \begin{array}{ll}
 \frac{\rho_{ji}+\rho_{ij}}{\sum_{\ell\in\cN_i}(\rho_{\ell i}+\rho_{i\ell})+\omega_i} = \frac{\rho_{ji}+\rho_{ij}}{\beta_{i}}, \!\!\!&\mbox{if}~e_{ij}\in \mathcal{A},\\
 \frac{\omega_i}{\sum_{\ell\in\cN_i}(\rho_{\ell i}+\rho_{i\ell})+\omega_i}= \frac{w_i}{\beta_{i}}, &\forall~i=j, \; i\in\cV,  \\
 0, &~\mbox{otherwise},
 \end{array}\right.\\
 {\Upsilon}& := \mbox{diag}\{{\beta}_1,\cdots, {\beta}_N\}\otimes I_M\succ 0.\label{eq:Upsilon}
\end{align}
Clearly $W$ is a row stochastic matrix satisfying the conditions in \eqref{eq:w:positive}. However, generally  $W$ constructed in this way is neither symmetric nor doubly stochastic, except when all $\beta_i$'s are identical. {We note that each entry of $W[i,j]$ is directly related to how agent $i$ will combine agent $j$'s information (this point will be made clear shortly).}

Surprisingly, Algorithm 2 admits a compact {\it single-variable} characterization, as we show in the following result.
\begin{proposition}\label{prop:equiv2}
{\it  The iteration \eqref{eq:x:update:static} -- \eqref{eq:lambda:update:static} of  Algorithm 2 (PGC) has the following compact characterization:
\begin{align}\label{eq:distributed}
&x^{r+1}-x^r +\Upsilon^{-1}(\zeta^{r+1}-\zeta^r)= \Upsilon^{-1}\left(-G\left(x^r\right)+G\left(x^{r-1}\right)\right)\nonumber\\
&+ (W\otimes I_{M}) x^r -\frac{1}{2}(I_{MN}+  W\otimes I_M)x^{r-1}, \;  \forall~r\ge 1,
\end{align}
where $\zeta^{r+1}\in\mathbb{R}^{MN}$ is a vector containing subgradients $\zeta^{r+1}_i\in\partial h_i(x_i^{r+1}),\;\forall~i\in \cV$.
In particular, each agent $i$ implements the following iteration: $\forall~r \ge 1$, {\small
\begin{align}\label{eq:distributed:i}
	&x^{r+1}_i-x^{r}_i + \frac{1}{\beta_i}(\zeta_i^{r+1}-\zeta^r_i)\! :=c^{r+1}_i \!=\! \frac{1}{\beta_i}\bigg(\!\!\!-\!\!\nabla g_i(x^r_i) \!+\! \nabla g_i(x^{r-1}_i)\!\!\bigg) \nonumber\\
	&\quad + \frac{1}{\sum_{j\in\cN_i}\hrho_{ij}+\omega_{i}/2}\bigg(\sum_{j\in\cN_i}\hrho_{ij} x^{r}_j +\frac{\omega_i}{2} x^r_i \bigg)\nonumber\\
	&\quad -  \frac{1}{2} \bigg(x^{r-1}_i + \frac{1}{\sum_{j\in\cN_i}\hrho_{ij}+\omega_i/2} \bigg(\sum_{j\in \cN_i} \hrho_{ij} x^{r-1}_j + \frac{\omega_i}{2} x^{r-1}_i\bigg) \bigg).
\end{align}}
}
\end{proposition}

The proof for the above claim is relegated to Appendix \ref{app:implementation}.  

{Let us take a closer look at iterations \eqref{eq:distributed} and \eqref{eq:distributed:i}. First, note that ${1}/{\beta_i}$ (or equivalently $\Upsilon^{-1}$) can be viewed as the {\it stepsize} for updating along the gradient direction. Second, for each node $i$, it is clear the penalty parameter $\hat{\rho}_{ij} = (\rho_{ij}+\rho_{ji})/2$ (or the $(i,j)$'s entry of the weight matrix $W$)  is the weight that specifies how the user $j$'s information (i.e., $x^r_j$ and $x^{r-1}_j$) is combined with the user $i$'s information at each iteration. The larger the value of $\hat{\rho}_{ij}$ (or $W[i,j]$),  the more emphasis that agent $i$ will put on agent $j$'s information.    
}

Then we comment on how \eqref{eq:distributed} and \eqref{eq:distributed:i} can be carried out in practice. 

If $h\equiv 0$, then $\zeta^r=0,\; \forall~r$. 
 To perform \eqref{eq:distributed:i} each agent $i$ needs its past iterate ($x^r_i$, $x^{r-1}_i$) the stepsize parameter $1/\beta_i$, the gradients ($\nabla g_i(x_i^r), \; \nabla g_i(x_i^{r-1})$), as well as the {\it weighted sum} of $x^r$ over its neighbors at the current and past iterations, i.e.,  $\sum_{j\in \cN_i} \hrho_{ij}x^r_j$ and $\sum_{j\in \cN_i} \hrho_{ij}x^{r-1}_j$, respectively. 
Also it is clear that the algorithm can be implemented in a fully distributed manner, since at iteration $r+1$, a given agent $i$ only communicates with its neighbors $\cN_i$.

 {When $h\ne 0$, iterations \eqref{eq:distributed} and \eqref{eq:distributed:i} can be implemented in the following manner. Assume that $x^0=x^{-1}=0$ and $\zeta^{0}=0$ for initialization. Then according to \eqref{eq:distributed:i} we have $x^1_i+\frac{1}{\beta_i}\zeta^1_i=0$, so $x^1_i$ and $\zeta^1_i$ can be obtained by solving the following problem
 $$x^1_i = \arg\min_{x_i}\; h_i(x_i)+\frac{\beta}{2} \|x_i\|^2 = \mbox{prox}^{\beta}_{h_i}(0).$$}
\!\!{Then one can compute $c^{2}_i$ according to \eqref{eq:distributed:i}}. To obtain $(x^{r+1}_i,\zeta^{r+1}_i)$,  $r\ge 1$, suppose $\zeta^r_i$ and {$c^{r+1}_i$ are}  available, 
then according to \eqref{eq:distributed:i},  we have
$$x^{r+1}_i+\frac{1}{\beta_i}\zeta^{r+1}_{i} = c^{r+1}_i+x^r_i+\frac{1}{\beta_i}\zeta^{r}_{i},$$
Finding $x^{r+1}_i$ is equivalent to solving the following
\begin{align}
x^{r+1}_i
&:={\rm prox}^{\beta_i}_{h_i}\left(c^{r+1}_i+x^r_i+\frac{1}{\beta_i}\zeta^{r}_{i}\right).\label{eq:prox}
\end{align}
Once $x^{r+1}_i$ is obtained, we can compute $\zeta^{r+1}_i$ by
\begin{align}
\zeta^{r+1}_i = \beta_i\left( x_i^r+c^{r+1}_i-x^{r+1}_i  \right) +\zeta^r_i.\label{eq:zeta}
\end{align}
Clearly, as long as problem \eqref{eq:prox} can be solved easily, iteration \eqref{eq:distributed} can be implemented efficiently in a distributed manner. 

\begin{table*}[t]
	\vspace{-0.3cm}
	\centering\small
	\caption{ Main Convergence Results.}\label{table:convergence}
	\begin{tabular}{*4c}
		\toprule
		\multicolumn{2}{c}{Scenario} & Convergence Condition & Convergence Rate \\
		\midrule
		Network Type & Gradient Type&{}  &{}  \\
		\midrule
		\midrule
		Static & Exact &  $\Omega+ \frac{1}{2}M_{+}(\Xi\otimes I_{M}) M^T_{+}-\widetilde{P}/2\succ 0$& $\mathcal{O}(1/r)$\\
		Static & Stochastic &  $\Omega+ \frac{1}{2}M_{+}(\Xi\otimes I_{M}) M^T_{+}-\widetilde{P}\succ 0$ & $\mathcal{O}(1/\sqrt{r})$\\
		Random & Exact  &  {$\Omega\succ \widetilde{P}/2$}  & $\mathcal{O}(1/r)$\\
		Random & Stochastic &  $\Omega\succ \widetilde{P}$ &$\mathcal{O}(1/\sqrt{r})$\\
		\bottomrule
	\end{tabular}
	\vspace{-0.4cm}
\end{table*}

\section{Convergence Analysis}
We begin analyzing the (rate of) convergence of the proposed methods. Let us define a diagonal matrix of Lipschitz constants by
\begin{align}\label{eq:tP}
\widetilde{P}:= \mbox{diag}\{P_1,\cdots, P_N\}\otimes I_{M}\in\mathbb{R}^{MN\times MN}.
\end{align}
Let $w:=[x; z; \lambda]$ denote the vector of primal-dual iterates generated by PGC/DySPGC, and let $w^*:=[x^*; z^*; \lambda^*]$ denote a vector of optimal primal-dual solutions for problem \eqref{eq:consensus:admm:matrix}.
Our main convergence results are summarized in Table \ref{table:convergence}. All the proofs of this section are relegated to the Appendix. 

\subsection{Analysis for Static Graphs}
For the PGC algorithm which use static graph and exact gradients, we have the following  convergence result.
\begin{theorem}\label{tm:time:invariant}
{\it Suppose that Assumption 1 holds,  $\cG^{r}=\cG$ for all $r$, and $\cG$ is connected.
Then the following hold:

\noindent {\bf \mbox{(a)}}\;  Algorithm 2 converges to a primal-dual optimal solution of \eqref{eq:consensus:admm:matrix} if the following condition is satisfied
\begin{align}
2\Omega+ M_{+}(\Xi \otimes I_{M}) M^T_{+}=\Upsilon {(W\otimes I_M)}+\Upsilon \succ \widetilde{P}\label{eq:graph:cond}.
\end{align}

\noindent {\bf \mbox{(b)}}\; Assume that $\mbox{dom}(h)$ is a bounded set, i.e., there exists a finite $C>0$ such that
$$d_{x}:=\sup_{\hat{x}, \; \tilde{x}\in\mbox{dom}(h)}\|\hat{x}-\tilde{x}\|\le C.$$
Suppose that $w^r:=[x^r; z^r; \lambda^r]$ is generated by Algorithm 2 and the stepsize matrix satisfies
\begin{align}
2\Omega+ M_{+}(\Xi\otimes I_{M}) M^T_{+}=\Upsilon { (W\otimes I_M)}+\Upsilon \succ 2\widetilde{P}.\label{eq:Omega:inexact}
\end{align}
Moreover, define
\begin{align*}
\barw^{r+1}&=[\bar x^r; \bar z^r; \bar \lambda^r]:=\frac{1}{r+1}\sum_{t=0}^{r} w^t, \\
d_z &:=
\hspace{-0.0cm}\sup_{\hat{x},\; \tilde{x}\in\mbox{dom}(h)}\sqrt{\sum_{ij:e_{ij}\in
\cA}2\rho_{ij}\|\hat{x}_i-\tilde{x}_j\|^2},
\end{align*}
and $d_{\lambda}(\rho):=\sup_{\lambda\in\mathcal{B}_\rho}\|\lambda-\lambda^0\|_{\Gamma^{-1}}^2$ where $\mathcal{B}_{\rho}:=\{\lambda\mid \|\lambda\|\le \rho\}$, for any  $\rho>0$.
Then for all $r>0$, we have
\begin{align}
&P(\bar{x}^r,\bar{z}^r):=f(\bar{x}^r)-f({x}^*)+\rho\|A \bar{x}^r+ B\bar{z}^r\| \notag \\
&~~~~~~~~~~~~\le \frac{1}{2r}\left(d^2_z+ d^2_{\lambda}(\rho)+\max_{i}\omega_i d^2_x\right). \label{eq:pgc:sub}
\end{align}
}
\end{theorem}

{Let us briefly comment on the assumptions made in each part of the above statement. In part {\bf (a)}, the condition \eqref{eq:graph:cond} imposes  requirements on the parameters of the algorithm, such as the proximal matrix $\Omega$ and the matrix $\Xi$ (which contains all penalty parameters $\{\rho_{ij}\}$). Note that a sufficient condition for \eqref{eq:graph:cond} is that $2\Omega\succ \widetilde{P}$, which is equivalent to $\omega_i>P_i/2$ for all $i\in\cV$.} 

{Also in part {\bf (b)}, $d_x$ represents the diameter of the feasible set $\mbox{dom}(h)$; $d_z$ can be viewed as the maximum size of any two $z$'s generated by the algorithm (see the first inequality in Appendix \ref{app:thm:time:variant:inexact:rate}); $d_\lambda(\rho)$ can be viewed as the distance between the initial solution $\lambda^0$ to the ball $\mathcal{B}_{\rho}$. Additionally, the boundedness of the set $\mbox{dom}(h)$ can be achieved for example when $h$ is the indicator function of certain bounded convex set. }

{The key novelty, as well as the main challenge,  in the analysis of the proposed approach is a careful bounding of the proximal parameter $\omega_i$, which results in faster practical numerical convergence performance (to be shown in Section \ref{sec: simulations}).} Indeed, compared with the existing convergence results on proximal-based ADMM such as \cite{gao14} and \cite{Ouyang15}, our bound for $\omega_i$ is reduced by at least a half. More importantly, {no global information is needed at each agent to verify such condition, in contrast to \cite{Ling15decentralize,shi14extra}}.  It is also interesting to note that the condition  \eqref{eq:graph:cond}, which only guarantees convergence,  is indeed {\it weaker} than the condition \eqref{eq:Omega:inexact}, which guarantees the global sublinear convergence rate.

%
%


Next we analyze the algorithm for static graph and stochastic gradient (i.e., Algorithm 1 applied to a static graph). 
\begin{theorem}\label{thm:time:variant:inexact:rate}
{\it Suppose that Assumption 1 holds, and the graph is static and connected (with $\cG^r=\cG$ for all $r$).
Suppose that $w^r$ is generated by Algorithm 1, and all the assumptions made Theorem \ref{tm:time:invariant}(b) hold true. If additionally the penalty parameter sequence $\{\eta^r\}$ satisfies
$$\eta^{r+1} = \sqrt{r+1},\quad \forall~r,$$
then it holds that
\begin{align*}
&\mathbb{E}\left[ P(\bar{x}^r,\bar{z}^r) \right]\le \frac{\sigma^2}{\sqrt{r}}+\frac{d^2_{x}}{2\sqrt{r}}+\frac{1}{2r}\left(d^2_z+ d^2_{\lambda}(\rho)+\max_{i}\omega_i d^2_x\right).
\end{align*}

}
\end{theorem}

\begin{remark}
In the previous two results, we have used $P(\bar{x}^r,\bar{z}^r)$
to measure the quality of the solution. This is a reasonable measure: according to \cite[Lemma 2.4]{gao14}, {when $\rho$ is large enough (in the sense that $\rho>\|\lambda^*\|$)}, $P(\bar{x}^r,\bar{z}^r)\le \epsilon$ implies that
$$ |f(\bar{x}^r)-f({x}^*)|\le \mathcal{O}(\epsilon), \quad \|A \bar{x}^r+ B\bar{z}^r\|\le \mathcal{O}(\epsilon).   $$
That is, both the constraint violation and the objective gap are in the same order as $\epsilon$. \hfill $\blacksquare$
\end{remark}


\begin{remark}\label{rmk:stochastic:admm}
	We remark that the stochastic ADMM method for solving general linearly constrained problem has been discussed in several recent papers \cite{Zhong14,wang13icassp,ouyang13,gao14}. However its application and the rate analysis in the context of distributed consensus based optimization appears to be new. {In particular, compared with the SGADM proposed in \cite{gao14}, our scheme only linearizes the objective function $f_i$, but not the entire augmented Lagrangian. Further, the order of the updates of the two primal variables has been reversed. These key differences make the analysis in \cite{gao14} not directly applicable.}
	\hfill $\blacksquare$
\end{remark}

\subsection{Analysis for Random Graphs}
In this section we analyze the convergence properties of Algorithm 1 (DySPGC) for random graphs defined in Definition \ref{def:graph}. The convergence claims are similar to those given in the previous section, but in the sense of convergence in expectation or with probability 1 (w.p.1).

We first analyze the simple case with exact gradient. 
To proceed, define a new function $J(x,z,\lambda)$ as
\begin{align}
J(x,z,\lambda) := \sum_{i=1}^{N}\frac{1}{\alpha_i} f_i(x_i)+\langle \lambda, A \Psi^{-1} x + B\Phi^{-1} z\rangle. \label{eq:J}
\end{align}
Define the following quantities
	\begin{align*}
	\tilde{d}_{x}: &= \min_{\hat{x},\tilde{x}\in\mbox{dom}(h)}(\tilde{x}-\hat{x})^T \Psi^{-1} (\tilde{x}-\hat{x}),\\
	\tilde{d}_{\lambda}(\rho): &= \sup_{\rho\in\mathcal{B}_{\rho}}(\lambda^0-\lambda)^T \Phi^{-1/2}\Gamma^{-1}\Phi^{-1/2}(\lambda^0-\lambda),
	\end{align*}
where $\Psi$ and $\Phi$ are given in \eqref{eq:Psi}. {These quantities can be viewed as matrix scaled versions of their counterparts \{$d_x, d_{\lambda}(\rho)$\} in the statement of Theorem \ref{tm:time:invariant}. }

The derivation of the following result is mostly based on that of Theorem \ref{tm:time:invariant}; the details can be found in our technical report \cite{hong15consensus}.

\begin{theorem}\label{thm:random:exact}
{\it Suppose that Assumption 1 holds, and  $\widetilde{G}(x^r, \xi^{r+1})=G(x^r), \; \forall~r$. Suppose that the graph $\{\cG^r\}$ is generated according to Definition \ref{def:graph}. 
	Then the following two statements hold true.

\noindent {\bf (a)} If the following condition is satisfied
\begin{align}
 2\Omega\succ \widetilde{P},\label{eq:random:stepsize}
 \end{align}
then $w^r$ generated by Algorithm 1 converges w.p.1. to a primal-dual solution of problem \eqref{eq:consensus:admm:matrix}.

\noindent {\bf (b)} Define $\bar{w}^{r}$ similarly as in the statement of Theorem \ref{tm:time:invariant}.
Suppose the following holds true
\begin{align}
\Omega\succ \widetilde{P}, \label{eq:random:stepsize:rate}
\end{align}
then Algorithm 1 generates a sequence $\bar{w}^r$ that satisfies
\begin{align*}
&\mathbb{E}\left[ P(\bar{x}^r,\bar{z}^r) \right]\le \frac{1}{2r}\left(2d_J+{d}^2_z+ \tilde{d}^2_{\lambda}(\rho)+\max_{i}\omega_i \tilde{d}^2_x\right),
\end{align*}
where $d_{J}:=\sup_{\lambda\in\mathcal{B}_\rho} J(x^0, z^0, \lambda) $.}
\end{theorem}
{Let us briefly compare the assumptions made in each of the statement. The condition $2\Omega\succ\widetilde{P}$ is equivalent to the condition that $\omega_i> P_i/2, \; \forall~i$, which implies that each local agent's proximal parameter should be chosen larger than $P_i/2$. Again, this condition is more relaxed compared with the one given in part (b), which requires a set of larger local proximal parameters $\{\omega_i\}$.}

It is interesting to note that the stepsize rules \eqref{eq:random:stepsize} and \eqref{eq:random:stepsize:rate} are both implied by their respective counterparts \eqref{eq:graph:cond} and \eqref{eq:Omega:inexact}, but the new rules are no longer related to the network structure.
Finally we analyze the case where the gradients are stochastic [i.e., Algorithm 1 (DySPGC) in its most general form]. 

\begin{theorem}\label{thm:random:inexact:rate}
{\it Define $\bar{w}^{r}$ similarly as in the statement of Theorem \ref{tm:time:invariant}.
Suppose that
$$\eta^{r+1} = \sqrt{r+1},\; \forall~r, \; \mbox{and}\quad  \Omega\succ \widetilde{P}.$$
Then Algorithm 1 generates a sequence $\bar{w}^r$ that satisfies
\begin{align*}
\mathbb{E}\left[ P(\bar{x}^r,\bar{z}^r) \right]\le &\frac{\sigma^2}{\sqrt{r}}+\frac{\max_{i}\omega_i(d^2_{x}+2\tilde{d}^2_x)}{2\sqrt{r}} \notag \\
&+\frac{1}{2r}\left(2d_J+{d}^2_z+ \tilde{d}^2_{\lambda}(\rho)+\max_{i}\omega_i d^2_x\right),
\end{align*}
where $d_{J}$ is defined in Theorem \ref{thm:random:exact}(b).}
\end{theorem}

The detailed proof can be found in our technical report \cite{hong15consensus}. {We note that comparing with the existing analysis for random graphs in \cite{Wei13,chang14}, our proof further takes into account inexact gradient information, and it does not assume any strong convexity on the objective functions. }

\section{Comparison with Existing Algorithms}\label{sec:compare}

Our proposed DySPGC as well as its special case PGC is closely related to a few existing algorithms. In this section we provide a detailed account of such relations; see Table \ref{table:compare} for a summary.

\vspace{-0.2cm}
\begin{table}[!htbp]
\centering\small
\caption{ Comparison of Different Algorithms with DySPGC.}\label{table:compare}
\centering
\begin{tabular}{lll}
\toprule
{Algorithm} & Connection to  & Special Setting \\
& DySPGC & \\

\midrule
IC-ADMM& Special Case  & Static, $\tG=\nabla g$, $g$ composite\\
DLM& Special Case  & Static, $h\equiv 0$, $W=W^T$, $\tG=\nabla g$\\
EXTRA& Special Case & Static, $h\equiv 0$, $W=W^T$, $\tG=\nabla g$\\
PG-EXTRA& Special Case  & Static, $W=W^T$, $\tG=\nabla g$\\
\midrule
DSG& Different $x$-step & N/A (not special case)\\
\bottomrule
\end{tabular}
\end{table}
\vspace{-0.2cm}

\subsection{Connection with the IC-ADMM}
Recently, an IC-ADMM algorithm is proposed in \cite{chang14distributed}, which solves the following problem in a distributed manner
\begin{align}\label{eq:consensus:2}
\min_{y\in\mathbb{R}^{M}}\; \sum_{i=1}^{N} \ell_i(A_iy) +h_i(y),
\end{align}
where $\ell_i(\cdot)$ is a \emph{strongly convex} function, and each $A_i$ is a given matrix not necessarily having full column rank. Clearly, this problem is a special case of our consensus problem \eqref{eq:consensus}, with the additional requirement that the smooth part of the objective has the composite form (strongly convex plus a linear mapping) given in \eqref{eq:consensus:2}.  The IC-ADMM algorithm is a special case of our Algorithm 2 (PGC) applied to solve problem \eqref{eq:consensus:2}, with constant penalty parameter $\rho_{ij}=\rho>0$, for all $i,j$. The analysis provided in \cite[Theorem 1]{chang14distributed} requires that the stepsize $1/\beta_i$ to be proportional to the strong convexity constant of the function $\ell_i(\cdot)$, which can be tiny for badly scaled functions. In our analysis, no such condition is necessary.

\vspace{-0.2cm}
\subsection{Connection with the DLM algorithm}
The Decentralized Linearized Alternating Direction Method of Multipliers (DLM) proposed in \cite{Ling15decentralize} is closely related to IC-ADMM. The DLM solves \eqref{eq:consensus} with $h_i\equiv 0$.
Its basic iteration is again Algorithm 2 (PGC) with parameters $\rho_{ij}=\rho>0$ and $\omega_i=\omega\ge 0$ for all $i,j$. The convergence condition in \cite[Theorem 1]{Ling15decentralize} is given by (described using our notation)
$$\beta \lambda_{\min}\big(\frac{1}{2}M_{+} M^T_{+}\big)+\omega > \max_i{P_i}/2.$$
This condition is an immediate consequence of the condition \eqref{eq:graph:cond} (with uniform $\rho_{ij}$'s and uniform $\omega_i$'s). 
%

\subsection{Connection with EXTRA}\label{sub:extra}
We show that Algorithm 2 (PGC) can be viewed as a generalization of the EXTRA \cite{shi14extra}. Consider applying Algorithm 2 (PGC) to problem \eqref{eq:consensus:admm:matrix} with a smooth objective (i.e., $h_i\equiv 0$ for all $i$). According to Proposition \ref{prop:equiv2}, one can write the iterates of 
Algorithm 2 as
\begin{align}
\!\!\!\!\!\!x^{r+1} = x^r + \Upsilon^{-1}\left(G(x^{r-1}) - G(x^{r})\right) + \widehat{W} x^r -\widetilde{W} x^{r-1},\label{eq:extra:form}
\end{align}
where
\begin{align} 
  \widehat{W}:= W\otimes I_M, \quad \widetilde{W}:= \frac{1}{2} (I_{MN}+W\otimes I_M). \label{eq:tW}
  \end{align}
{Eq. \eqref{eq:extra:form} is precisely the EXTRA update developed in \cite{shi14extra}, except for the two relatively minor points:
\begin{enumerate}
\item In \eqref{eq:extra:form} a slightly more general {\it matrix stepsize}  $\Upsilon^{-1}$ is used instead of the scalar stepsize used in EXTRA.
\item The EXTRA allows a slightly wider choice of $ \widetilde{W}$, i.e., $ {1}/{2} (I_{MN}+W\otimes I_M)\succeq \widetilde{W}\succeq W$, $\mbox{null}\{W-\widetilde{W}\} =\mbox{span}\{\mathbf{1}\}$ and $\mbox{null}\{I_{MN}-\widetilde{W}\} \supseteq \mbox{span}\{\mathbf{1}\}$, where $\mathbf{1}$ is an all one vector of appropriate size.  However except for the common choice  \eqref{eq:tW}, these conditions are difficult (if not impossible) to verify in a fully distributed manner.
\end{enumerate}}
When a single scalar stepsize is used (as was done in EXTRA), say $\beta=\beta_i=\beta_j>0$ for all $i,j$, then we can perform either one of the following procedures to identify the parameters of Algorithm 2 (PGC) (depending on whether the weight matrix $W$ is known {\it a priori}):

\noindent{\bf From Algorithm Parameters to Weight Matrix.} Suppose the agents can select $\{\omega_i\}$ and $\{\rho_{ij}\}$. Then for any set of fixed $\{\rho_{ij}\}$'s, pick $\beta$ and $\omega_i$'s such that
$$ \omega_i = \beta - \sum_{j\in\cN_j} (\rho_{ij}+\rho_{ji})\ge0, \quad \forall~i,$$
(c.f.  \eqref{eq:beta}). Further, pick $\beta$ large enough such that convergence conditions such as \eqref{eq:graph:cond} are satisfied.  Note that the weight matrix in \eqref{eq:W} induced by such choice of parameters must be symmetric and doubly stochastic.\\
\noindent {\bf From Weight Matrix to Algorithm Parameters.} Suppose the weight matrix $W$ is given and fixed, and it is a symmetric doubly stochastic matrix. The symmetry of $W$ implies $\beta_i=\beta_j=\beta$. For any fixed $\beta>0$, one can easily find the parameters $\{\rho_{ij}\}$ and $\{\omega_i\}$ by letting $\rho_{ij}+\rho_{ji}=\beta\times W[i,j]$, $\omega_i = \beta\times W[i,i]$ for all $ij$ such that $e_{ij}\in\cA$. Again one should pick $\beta$ large enough such that the convergence conditions (i.e., \eqref{eq:graph:cond}) are satisfied. Note that such construction implies that $\sum_{\ell\in\cN_i}(\rho_{\ell i}+\rho_{i\ell})+\omega_i= \sum_{\ell\in\cN_i}\beta W[i, \ell]+\beta W[i,i]=\beta$ (since $\sum_{j\in\cN_i}W[i,j] + W[i,i]=1$), which recovers its original definition in \eqref{eq:beta}.

To compare the convergence result in Theorem \ref{tm:time:invariant} and that of \cite[Theorem 3.3]{shi14extra}, note that when the scalar stepsize is used, we have $\Upsilon = \beta I_{MN}$. Therefore a sufficient condition to guarantee the condition  in Theorem \ref{tm:time:invariant} is that
\begin{align*}
\beta\lambda_{\min}\left(I_{MN}+W\otimes I_M\right)> \max_{i} P_i.
\end{align*}
This is precisely the condition set forth in \cite[Theorem 3.3]{shi14extra}.

From the above expression it is clear that  $\beta$  depends on \emph{all} the local functions, therefore it has to be decided in a centralized manner. In contrast, the stepsize parameters in  PGC can be chosen as: $\omega_i\ge P_i/2$ (cf. the remarks made after Theorem \ref{tm:time:invariant}). The latter choice is simple, distributed implementable, and more importantly it results in improved convergence speed in practice, especially when the curvatures of $g_i$'s vary significantly, i.e., $\max_i{P_i}\gg\min_i{P_i}$. This will be demonstrated in Section \ref{sec: simulations}. 

{We comment that in a couple of recent works \cite{Mokhtari16linear} and \cite{Mokhtari16DSA}, the authors have established that EXTRA is also related (and in fact in most cases equivalent) to certain saddle point method, and certain proximal augmented Lagrangian method. Combining the observation made in this work, we can conclude that all these methods (i.e., the saddle point method \cite{Mokhtari16DSA}, the proximal augemented Lagrangian method \cite{Mokhtari16linear}, the EXTRA and the PGC) are all closely connected \footnote{We thank the anonymous reviewer for bringing these new developments to our attention.}. }




\vspace{-0.2cm}
\subsection{Connection with PG-EXTRA}
One can also show that the proposed Algorithm 2 (PGC) generalizes the PG-EXTRA \cite{shi15pgextra}.
According to the argument leading to \eqref{eq:prox}, one can explicitly express \eqref{eq:distributed} by 
\begin{align}
	x^{r+1}_i&\stackrel{\eqref{eq:prox}}={\rm prox}^{\beta_i}_{h_i}\big(c^{r+1}_i+x^r_i+\frac{1}{\beta_i}\zeta^{r}_{i}\big)\nonumber\\
	&\stackrel{\eqref{eq:zeta}} = {\rm prox}^{\beta_i}_{h_i}\big(c^{r+1}_i+ c^{r}_i+ x^{r-1}_i+\frac{1}{\beta_i}\zeta^{r-1}_{i}\big)\nonumber\\
	&= {\rm prox}^{\beta_i}_{h_i}\bigg(\sum_{t=2}^{r+1}c^{t}_i+ x^1_i+\frac{1}{\beta_i} \xi^1_i\bigg)\nonumber.
	\end{align}
By the definition of $c_i$ in \eqref{eq:distributed:i}, we have  

	\begin{align*}
	\sum_{t=2}^{r+1}c^{t}_i&=\frac{1}{\beta_i}\left(-\nabla g_i(x^r_i)+\nabla g_i(x^0_i)\right)+ \widehat{W}_i x^r\nonumber\\
	&\quad +\sum_{t=2}^{r}(\widehat{W}_i-\widetilde{W}_i)x^{t-1}-\widetilde{W}_i x^0
	\end{align*}
where $\widehat{W}_i$ and $\widetilde{W}_i$ denote the $i$th row of $\widehat{W}$ and $\widetilde{W}$ (as have been defined in \eqref{eq:tW}), respectively. Again by \eqref{eq:distributed:i}, and assume that $x^0=0$ and $\nabla g_i(x^{-1}_i)=0$, we can check that
\begin{align*}
& x^1_i+\frac{1}{\beta_i} \xi^1_i=-\frac{1}{\beta_i}\nabla g_i(x_i^0).
\end{align*}
Combining the above three equalities we have
	\begin{align}
	\hspace{-0.2cm}x^{r+1}_i\hspace{-0.2cm} = {\rm prox}^{\beta_i}_{h_i}\bigg(\frac{-1}{\beta_i}\nabla g_i(x^r_i)+ \widehat{W}_i x^r+\sum_{t=1}^{r}(\widehat{W}_i-\widetilde{W}_i)x^{t-1}\bigg). \notag 
	\end{align} This is the PG-EXTRA proposed in \cite[Algorithm 1]{shi15pgextra}.
%

\vspace{-0.2cm}
\subsection{Connection with the DSG Method}\label{sub:distributed:gradient}
Below we show that Algorithm 2 (PGC) is closely related to the DSG iteration \eqref{eq:distributed:subgradient:method}.
{Assume for simplicity that $h_i\equiv 0$ for all $i$.} Suppose that the $z$ and $\lambda$ steps of the PGC remain the same while the $x$-step \eqref{eq:x:update:static} is replaced by the following

\begin{align*}
x^{r+1} \hspace{-0.0cm}=\arg\min_{x} \bigg\{ \left\langle G(x^r), x-x^r\right\rangle &+ \frac{1}{2}\left\|Ax+Bz^r\right\|_{\Gamma}^2 \notag \\
&+\frac{1}{2}\|x-x^r\|^2_{\Omega}\bigg\}.
\end{align*}
That is, in the $x$-step we let $\lambda^r= 0$. The claim is that by such modification one recovers the DSG iteration \eqref{eq:distributed:subgradient:method}.
To argue this, we write down the optimality condition of the modified iteration as
\begin{subequations}
\begin{align*}
&G(x^r) + A^T\Gamma (A x^{r+1}+ B z^r) + \Omega (x^{r+1}-x^r) = 0,\\
& B^T \lambda^r +  B^T \Gamma (A x^{r+1}+ B z^{r+1}) = 0, \\
& \lambda^{r+1}-\lambda^r -\Gamma  \left( A x^{r+1}+ B z^{r+1}\right) = 0 .
\end{align*}
\end{subequations}
Following the derivation of Proposition \ref{prop:equiv2} until \eqref{eq:opt:x:2} in Appendix \ref{app:implementation}, we have
\begin{align}
&G(x^r)+ \alpha^{r+1}-\alpha^{r} + \frac{1}{2}M_{+}(\Xi\otimes I_M) M^T_{+} (x^{r+1}-x^r)  \nonumber\\
&\quad\quad\quad + \Omega  (x^{r+1}-x^r) = 0, \nonumber\\
&\alpha^{r+1} = \alpha^{r}+\frac{1}{2}M_{-}(\Xi\otimes I_M) M^T_{-} x^{r+1}. \nonumber
\end{align}
Note that compared with \eqref{eq:opt:x:2}, the first equality above has an additional term $-\alpha^r$. Plugging the second equality  into the first one, we obtain
\begin{align*}
&G(x^r)+ \frac{1}{2}M_{-}(\Xi\otimes I_M) M^T_{-} x^{r+1} \nonumber\\
&\quad+ \frac{1}{2}M_{+}(\Xi\otimes I_M) M^T_{+} (x^{r+1}-x^r) + \Omega  (x^{r+1}-x^r) = 0.
\end{align*}

{By the definition of the matrices $M_{+}$ and $M_{-}$ in \eqref{eq:def:M}, one can verify the following identities 
{\small
	\begin{subequations}\label{eq:M:Xi:product}
		\begin{align}
		M_{+}(\Xi\otimes I_M) z & = \left[\begin{array}{l}
		\sum_{j\in \cN_1} (\rho_{j1} z_{j1}+\rho_{1j}z_{1j}) \\
		\sum_{j\in \cN_2} (\rho_{j2} z_{j2}+\rho_{2j}z_{2j}) \\
		\quad\quad\vdots\\
		\sum_{j\in \cN_N} (\rho_{jN} z_{jN}+\rho_{Nj}z_{Nj})
		\end{array}\right],\\
		M_{+}(\Xi\otimes I_M) M^T_{+} x & = \left[\begin{array}{l}
		\sum_{j\in \cN_1} (2\hrho_{j1} x_1 +2\hrho_{j1} x_j)\\
		\sum_{j\in \cN_2} (2\hrho_{j2} x_2 +2\hrho_{j2} x_j)\\
		\quad\quad{\vdots}\\
		\sum_{j\in \cN_N} (2\hrho_{jN} x_N +2\hrho_{jN} x_j)
		\end{array}\right], \label{eq:M:Xi:product:b}
		\end{align}
		\begin{align}
		M_{-}(\Xi\otimes I_M) M^T_{-} x & = \left[\begin{array}{l}
		\sum_{j\in \cN_1} (2\hrho_{j1} x_1 -2\hrho_{j1} x_j)\\
		\sum_{j\in \cN_2} (2\hrho_{j2} x_2 -2\hrho_{j2} x_j)\\
		\quad\quad\vdots\\
		\sum_{j\in \cN_N} (2\hrho_{jN} x_N -2\hrho_{jN} x_j)
		\end{array}\right]. \label{eq:M:Xi:product:c}
		\end{align}
	\end{subequations}}}
Utilizing \eqref{eq:M:Xi:product}, and by the definition of $\beta_i$ \eqref{eq:beta} 
and the definition of the weight matrix $W$ in \eqref{eq:W}, we can write the above iteration compactly as
\begin{align*}
x^{r+1} 
&= -\Upsilon^{-1}G(x^r)+\frac{1}{2}\left(I_{MN}+W\otimes I_M\right)x^r.
\end{align*}
After picking a uniform scalar stepsize $\beta_i=\beta_j=\beta>0$ (cf. Section \ref{sub:extra} for how this can be done), we immediately get the DSG iteration \eqref{eq:distributed:subgradient:method}  [with a weight matrix given by $\widetilde{W}=\frac{1}{2}(I_{MN}+W\otimes I_M)$].

Obviously, our convergence analysis does not work for this variant, as the $x$-update is no longer related to the dual variable $\lambda$. Indeed, to prove convergence of the DSG, an iteration-dependent and increasing $\beta$ is needed, and such convergence is usually slower than $\mathcal{O}(1/r)$; see \cite{shi14extra,Nedic09, chen12thesis,Jakovetic14} and the references therein. Nevertheless, the above observation reveals a fundamental connection between the ADMM-based method and the classical DSG method.

{
	\section{Extension to Accelerated DySPGC}
	The relationship identified between the DySPGC and the EXTRA, PG-EXTRA, IC-ADMM etc. provides a systematic way to analyze and generalize various existing algorithms. In this section, we provide one such generalization which accelerates the DySPGC (hence the EXTRA, PG-EXTRA, IC-ADMM, etc). The algorithm is inspired by  \cite{Ouyang15}.
	
	For simplicity, we will restrict ourselves to the static graphs  in this section. Let $\{\eta^{r}, \theta^{r}, \nu^{r}\ge 0\}$ denote a sequence of iteration-dependent parameters, whose values will be given shortly; Let $\{x^{r,\md}, x^{r,\ag}, z^{r,\ag}, \lambda^{r, \ag}\}$ denote a sequence of auxiliary variables. The proposed accelerated algorithm is given in the table below.

	\begin{center}
		\fbox{
			\begin{minipage}{5in}
				\smallskip
				\centerline{\bf Algorithm 3.  Accelerated DySPGC Over Static Graphs}\label{eq:acc}
				\smallskip
				At iteration $0$, let $B^T \lambda^0 = 0$, $z^0 = z^{0,\ag}=\frac{1}{2} M^T_{+} x^0$.
				
				At each iteration $r+1$, update the variable blocks by:
				\begin{subequations}
					\begin{align}
					x^{r+1,\md}& = (1-\nu^r)x^{r,\ag} +\nu^r x^r\\
					x^{r+1}& =\arg\min_{x}\; \bigg\{\left\langle \tG(x^{r+1,\md},\xi^{r+1}), x-x^r\right\rangle + h(x) \label{eq:x:update5}\\\
					&\hspace{-0.5cm}+ \frac{1}{2}\left\|Ax+Bz^r+\Gamma^{-1}{\lambda^r}\right\|_{\Gamma}^2+\frac{1}{2}\|x-x^r\|^2_{\theta^r\Omega+{\eta^{r+1}}I_{MN}} \bigg\}\nonumber\\
					x^{r+1,\ag}& = (1-\nu^r)x^{r,\ag} +\nu^r x^{r+1}\label{eq:amd:update5}\\
					z^{r+1} & =\arg\min_{z}\; \frac{1}{2}\left\|A x^{r+1} + Bz+\Gamma^{-1}{\lambda^r}\right\|_{\Gamma}^2\label{eq:z:update5}\\
					z^{r+1,\ag}& = (1-\nu^r)z^{r,\ag} +\nu^r z^{r+1}\label{eq:zag:update5}\\
					\lambda^{r+1}& = \lambda^r +\Gamma\left( A x^{r+1}+ Bz^{r+1}\right)\label{eq:lambda:update5}\\
					\lambda^{r+1,\ag}& = (1-\nu^r)\lambda^{r,\ag} +\nu^r \lambda^{r+1}\label{eq:lambdaag:update5}
					\end{align}
				\end{subequations}
				
			\end{minipage}
		}
	\end{center}

	First note that $x^{r,\ag},  z^{r,\ag}, \lambda^{r,\ag}$ are convex combinations of all previous iterates $\{x^{t}\}_{t=1}^{r}$, $\{z^{t}\}_{t=1}^{r}$, $\{\lambda^{t}\}_{t=1}^{r}$, respectively. Second, $x^{r+1,\md}$ is an intermediate point on which the stochastic gradient is evaluated. Therefore in total there are three sequences related to the $x$ update, resembling the Nesterov's acceleration scheme \cite{Nesterov04}.

	{The convergence rate of Algorithm 3 can be analyzed similarly as in \cite{Ouyang15}, we include the proof in the Appendix for completeness}. Compared with the bound given in Theorem \ref{thm:time:variant:inexact:rate}, the accelerated version is able to significantly reduce the scaling with respect to $\max_i w_i$,  which  in turn depends on the network structure as well as the Lipschitz constants of the local gradients through \eqref{eq:Omega:inexact}. 

	
	\begin{theorem}\label{thm:time:variant:inexact:rate:acc}
		{\it Suppose that the assumptions made in Theorem \ref{thm:time:variant:inexact:rate} are true. Further let
			\begin{align}
			\nu^{r} = \frac{2}{r+1}, \; \theta^r =\frac{2}{r+1}, \; \eta^r=\sqrt{r+1}, \;\varpi^{r} = \frac{2}{r(r+1)}.  \label{eq:nu}
			\end{align}
			Assume that the stepsize matrix satisfies
			\begin{align}
			4\Omega+ M_{+}(I_{M}\otimes\Xi) M^T_{+}\succ 4\widetilde{P}.\label{eq:Omega:inexact:2}
			\end{align}
			
			Then the iterates generated by Algorithm 3 satisfy{
				\begin{align*}
				&\mathbb{E}\big[f(x^{r+1,\ag})-f(x^*) +\rho\|Ax^{r+1,\ag}+Bz^{r+1,\ag}\|\big]\nonumber\\
				&\le \frac{1}{r}\big(d^2_z+d^2_{\lambda}(\rho)+ \frac{1}{r+1}\max_{i}\omega_i d^2_x\big)+ {\frac{2\sigma^2}{3}\frac{1}{\sqrt{r+1}}}+\frac{\sqrt{r+1}}{r}d^2_x.
				\end{align*}}
			where {$d_{\lambda}(\rho):=\sup_{\lambda\in\mathcal{B}_\rho}\sup_{\tilde{\lambda}}\|\lambda-\tilde{\lambda}\|_{\Gamma^{-1}}^2$, $\mathcal{B}_{\rho}=\{\lambda\mid \|\lambda\|\le \rho\}$}, and $\rho>0$ is any finite constant.
			
		}
	\end{theorem}
	
}

\begin{figure*}[!t]
			\vspace{-0.4cm}
	\begin{center}
		{\subfigure[][Static network, exact gradient]{\resizebox{.45\textwidth}{!}
				{\includegraphics{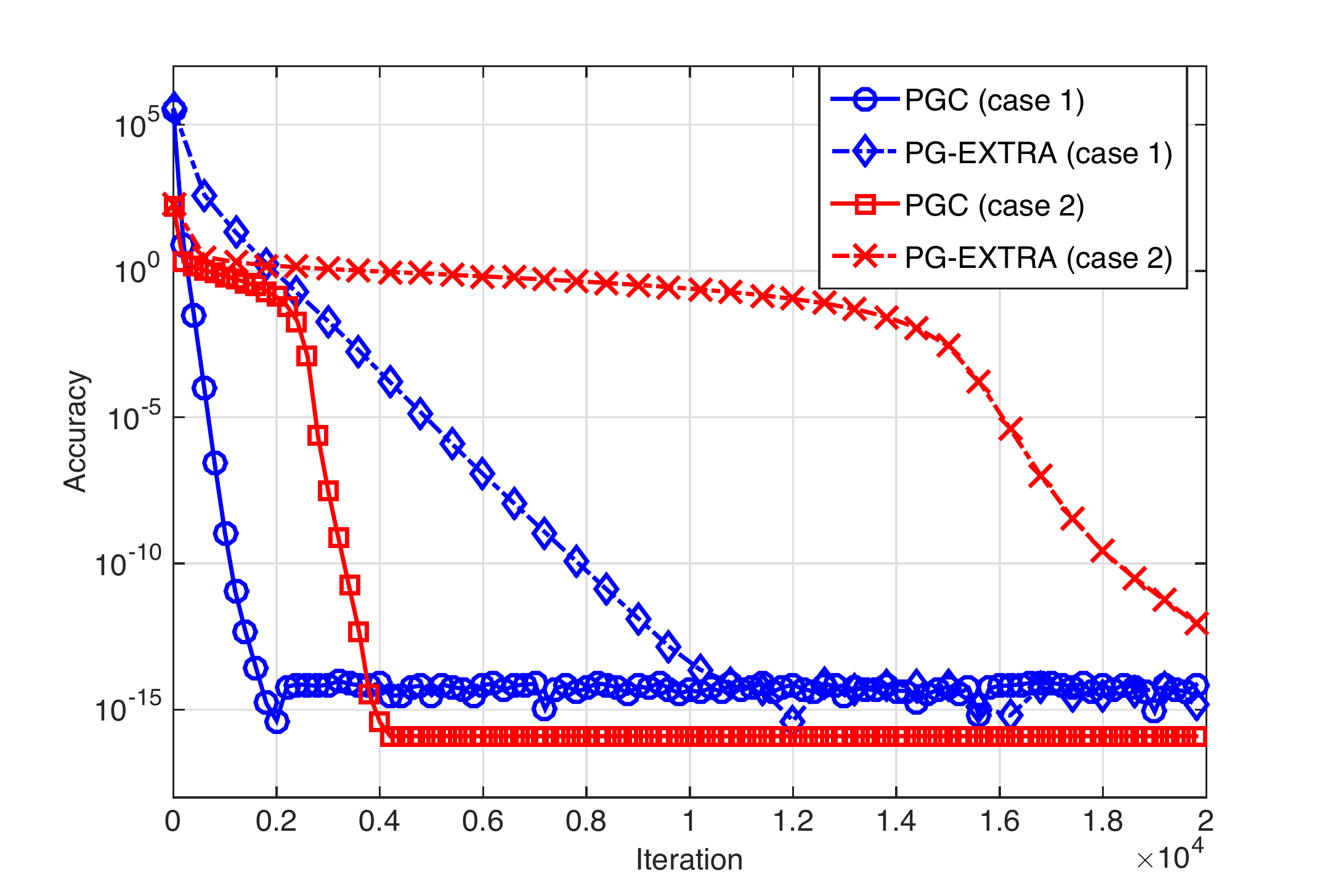}}} } \vspace{-0.1cm}
		{\subfigure[][Static network, exact gradient]{\resizebox{.45\textwidth}{!}{\includegraphics{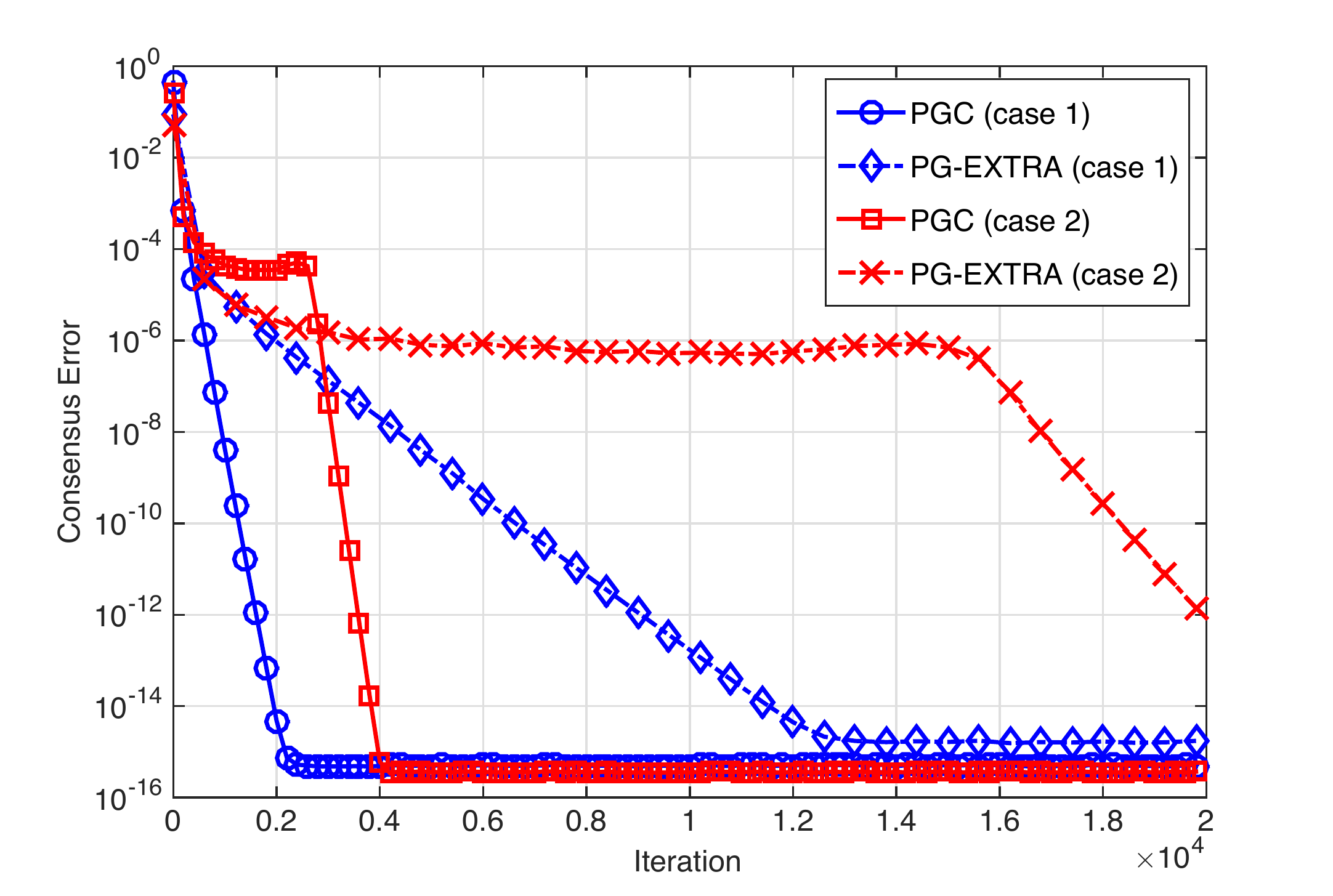}}}
		}
	\end{center}\vspace{-0.2cm}
	\caption{Convergence curves of proposed PGC algorithm  (Algorithm 2) with exact gradient information over a static network, for Case1:$M=1000,K=200,\nu=0.1$ and Case 2:  $M=1000,K=50,\nu=50$.		
	}
	\vspace{-0.2cm}\label{fig: spgc no noise}
\end{figure*}

\begin{figure*}[!h]
	\begin{center}
		\vspace{-0.2cm}
		{\subfigure[][Static network, exact gradient, $\nu=0$]{\resizebox{.45\textwidth}{!}
				{\includegraphics{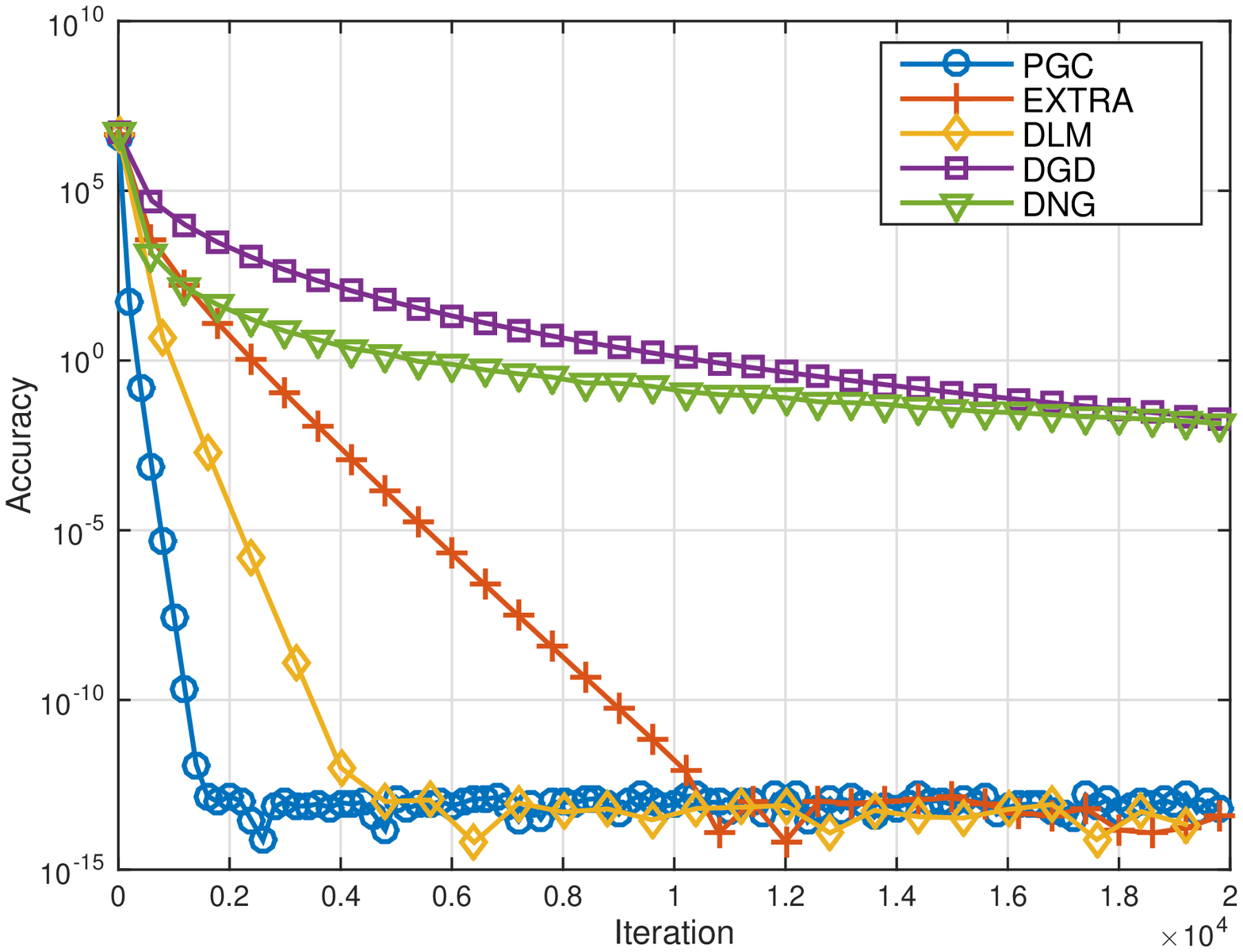}}} } \hspace{0.1cm}
		{\subfigure[][Static network, exact gradient, $\nu=0$]{\resizebox{.45\textwidth}{!}{\includegraphics{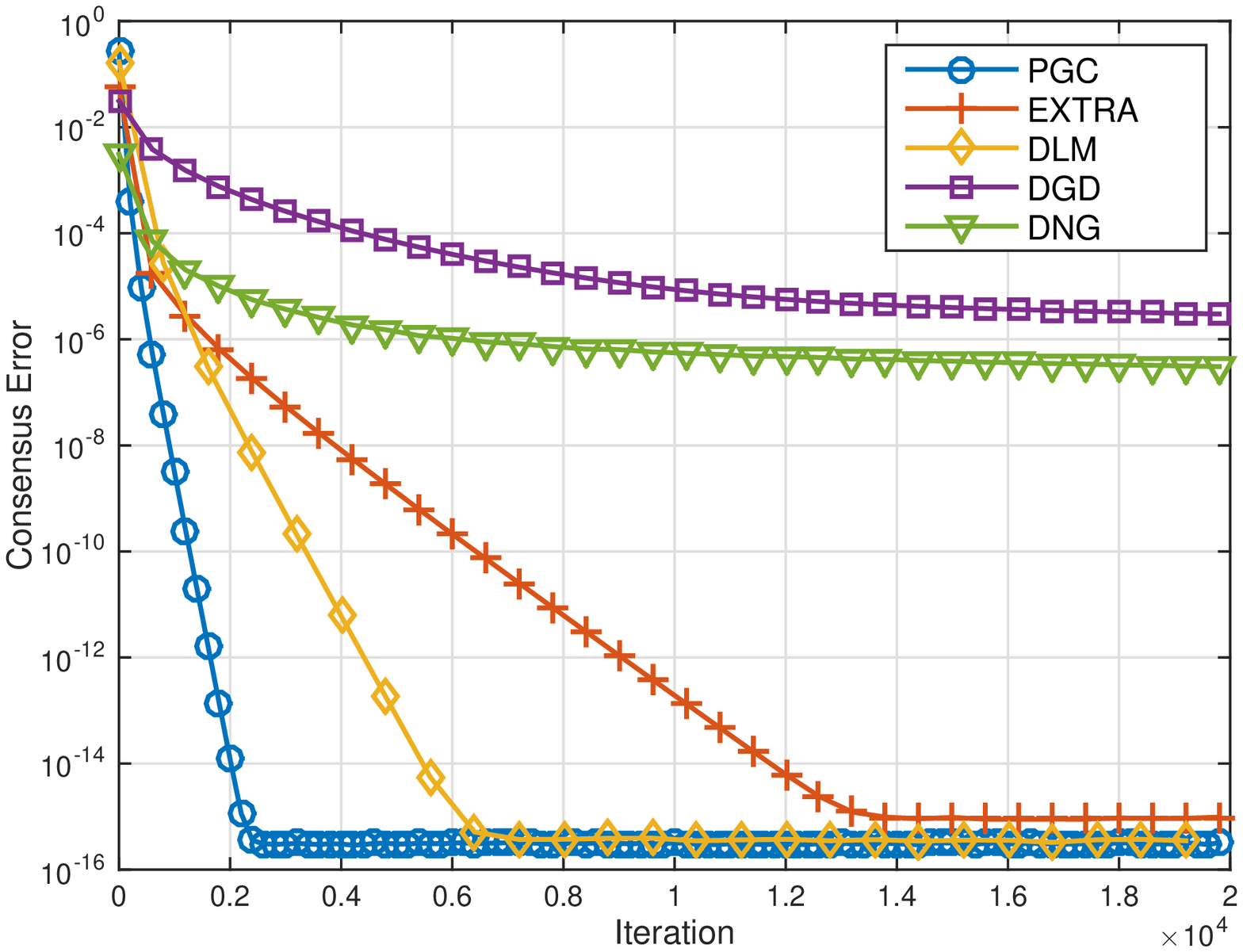}}}
		}
	\end{center}\vspace{-0.2cm}
	\caption{Convergence curves of proposed PGC algorithm  (Algorithm 2) with exact gradient information over a static network, for $M=1000,K=200,\nu=0$ (Least squares problem)  }
	\vspace{-0.2cm}\label{fig: pgc vs existing methods}
\end{figure*}

\vspace{-0.2cm}
\section{Numerical Results}\label{sec: simulations}
In this section, we present some simulation results of the proposed algorithms by solving the following LASSO problem
	\begin{align}\label{eqn: lasso}
	\min_x \textstyle \frac{1}{2}\sum_{i=1}^N  \|A_i x - b_i\|^2 + \nu\|x\|_1
	\end{align}
where $A_i\in\mathbb{R}^{K\times M}$, $b_i\in\mathbb{R}^{K}$ and $\nu>0$ is a penalty parameter. Each data matrix $A_i$ is randomly generated as $A_i=L_i\times Q_i$ where $L_i\sim\mbox{Uniform}[0,~10]$ and $Q_i\in\mathbb{R}^{K\times M}$ whose entries are i.i.d. standard Gaussian random variables; $b_i= A_i c + d_i$
where $c\in\mathbb{R}^M$ is a sparse random vector with $5$ percent of uniformly distributed non-zero entries; $d_i\in\mathbb{R}^K$ is a vector of i.i.d. zero-mean Gaussian random variables with standard deviation $0.01$. Note that here the Lipschitz constants are $P_i=\|A_iA^T_i\|,\; \forall~i$. The static graph $\mathcal{G}$ contains $16$ nodes  $|\cV|=N=16$, and the edge set $\cE$ is randomly generated following \cite{YildizScag08}, with a radius parameter set to $0.4$.

{We compare Algorithm 2 (PGC) with PG-EXTRA in \cite{shi15pgextra}, EXTRA in \cite{shi14extra}, DLM  in \cite{Ling15decentralize}, distributed gradient descent (DGD) algorithm in \cite{Nedic09subgradient}
and the distributed Nesterov gradient descent (DNG) algorithm in \cite{Jakovetic14}.
We also compare the static version of Algorithm 1 (i.e., SPGC) with a distributed stochastic gradient descent method (D-SGD) \cite{Ram2010_stoc}. The stepsize for the EXTRA/PG-EXTRA is chosen according to the sufficient condition suggested in \cite[Theorem 3.3]{shi14extra}, and the weight matrix $W$ is the Metropolis constant edge weight matrix. 
For Algorithm 1 (resp. Algorithm 2), $\omega_i={P_i}/{2}$ (resp. $\omega_i={P_i}$) and $\rho_{ij}=10^{3}$ for all $i,j$. 
For DLM, the parameter $c$ in \cite[Eqn. (21)]{Ling15decentralize} (which is equivalent to $\rho_{ij}$ here) is set to $10^3$, and $\rho$ (which corresponds to $\omega_i$ here) is set such that $\xi$ in \cite[Eqn. (21)]{Ling15decentralize} equals zero. 
For the DGD, DNG and D-SGD algorithms, the Metropolis weight matrix is used and 
the stepsize is set as $0.01/(r+5000)$
\footnote{The parameters are not searched in an exhaustive manner. Instead, for example, the parameter 
	 $\rho_{ij}$ in Algorithms 1 and 2 are tested for values 0.1, 1, 10, 50, 100, 500, 1000, 5000 and so on, and the value 1000 is chosen as the algorithm yields the fastest convergence performance among the others.}.  
}
To measure the progress of different algorithms, we define the following two quantities
\begin{align}
& \mbox{{Accuracy}} \textstyle = \frac{|f({x}^r) - f^*|}{f^*}, \notag \\ 
\vspace{-0.2cm}
&{\mbox{{Consensus Error}}  \textstyle=  \sqrt{\sum_{i=1}^N \|x_i^r - \hat{x}^r \|^2}/N },\nonumber
\end{align}
where $ f^*$ is the optimal objective value of problem \eqref{eqn: lasso} and is obtained by the FISTA method \cite{Beck:2009:FIS:1658360.1658364}, {and $\hat x^r = \frac{1}{N}\sum_{i=1}^N x_i^r $.}

\begin{figure*}[!t]
	\begin{center}		
		\vspace{-0.2cm}
		{\subfigure[][Static network, stochastic gradient]{\resizebox{.45\textwidth}{!}
				{\includegraphics{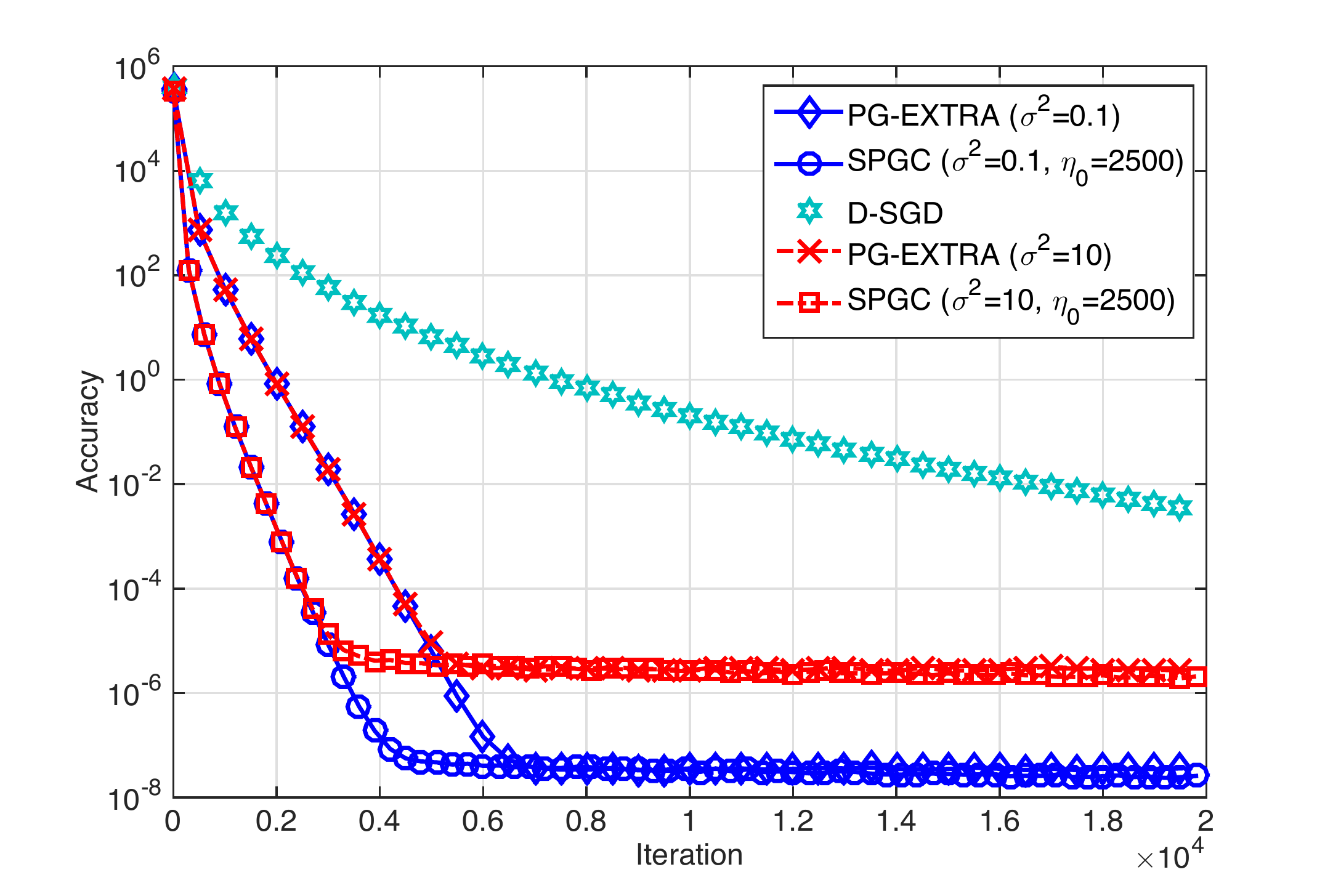}}} } \vspace{-0.1cm}
		{\subfigure[][Static network, stochastic gradient]{\resizebox{.45\textwidth}{!}{\includegraphics{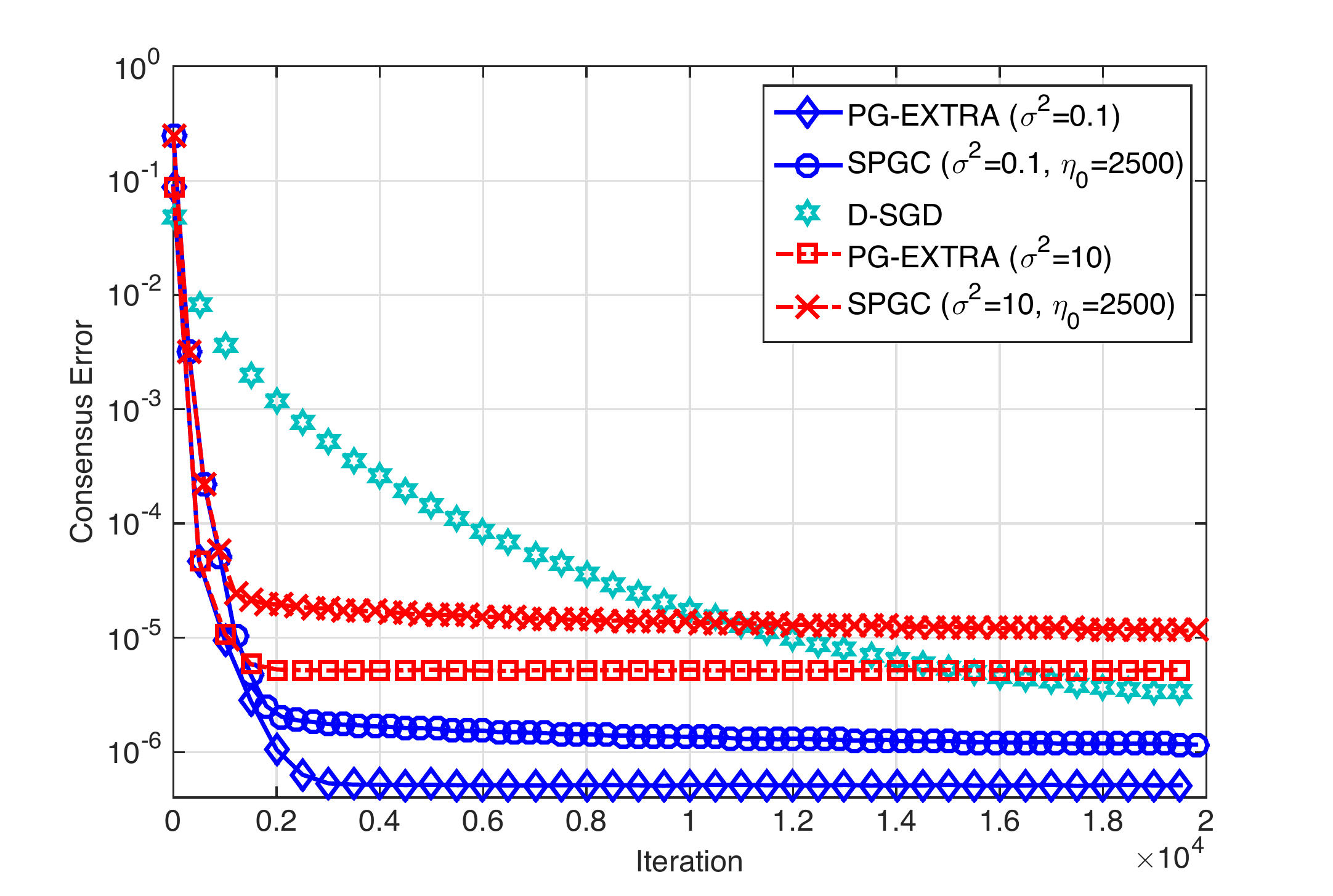}}}
		}
	\end{center}\vspace{-0.2cm}
	\caption{Convergence curves of proposed SPGC algorithm (Algorithm 1) with stochastic gradient information over a static network, for $M=1000,K=200,\nu=0.1$ and $\eta^r=\eta_0 \sqrt{r}$.  }
	\vspace{-0.4cm}\label{fig: spgc noise}
\end{figure*}

\begin{figure*}[!h]
	\begin{center}
		\vspace{-0.2cm}
		{\subfigure[][]{\resizebox{.45\textwidth}{!}
				{\includegraphics{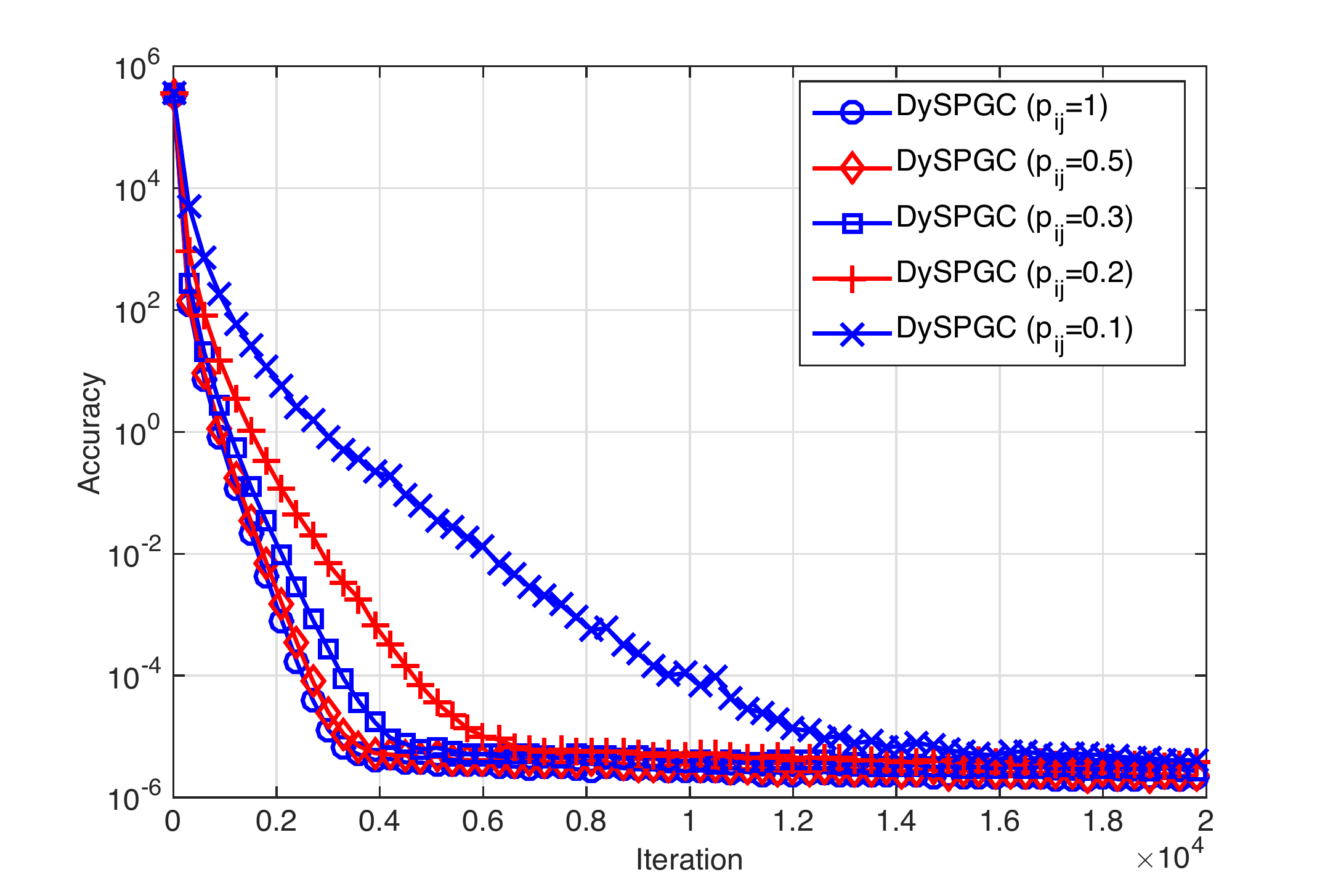}}} } \vspace{-0.1cm}
		{\subfigure[][]{\resizebox{.45\textwidth}{!}{\includegraphics{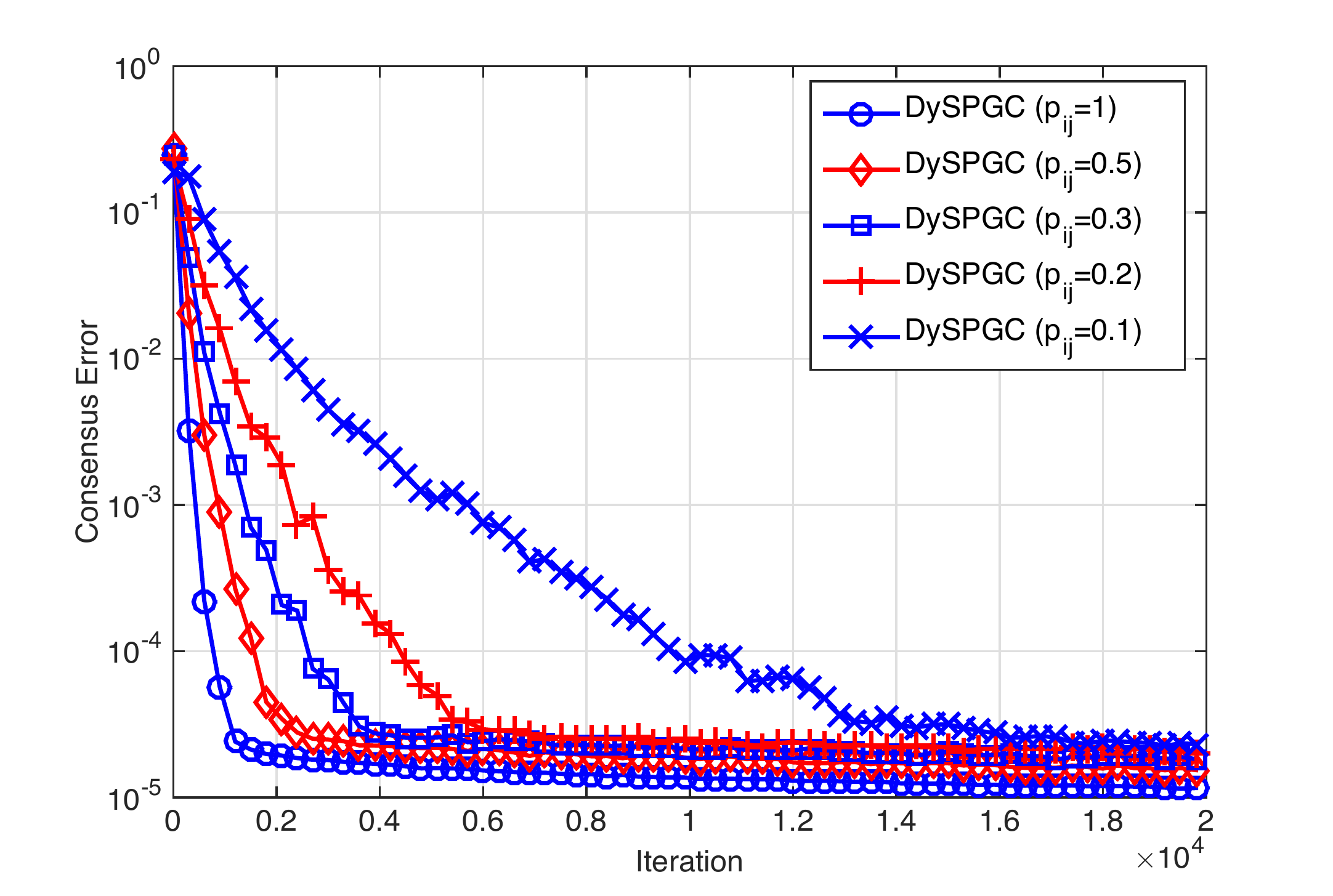}}}
		}
	\end{center}\vspace{-0.2cm}
	\caption{Convergence curves of proposed DySPGC algorithm with stochastic gradient information ($\sigma^2=0.1$) over a time-varying network, for $M=1000,K=200,\nu=0.1$ and $\eta^r=2500 \sqrt{r}$.   }
	\vspace{-0.2cm}\label{fig: spgc time varying}
\end{figure*}

The convergence curves of the proposed PGC (Algorithm 2) and the PG-EXTRA \cite{shi15pgextra} are shown in Figure \ref{fig: spgc no noise}, by assuming static networks and exact gradient information.  Two problem settings are considered: {\bf Case 1).} $M=1000,K=200,\nu=0.1$ and {\bf Case 2).}  $M=1000$, $K=50$, $\nu=50$. Given $N=16$, problem \eqref{eqn: lasso}  is strongly convex for Case 1 and (non-strongly) convex for Case 2. 
One can see from Figure \ref{fig: spgc no noise} the proposed algorithm outperforms the PG-EXTRA, in terms of both accuracy and consensus error. 
This is expected since compared with the PG-EXTRA, the PGC is able to use larger and more flexible stepsizes, as discussed in Section \ref{sub:extra}.

{As EXTRA, DLM, DGD and DNG are developed for smooth problems,  we set $\nu=0$, $M=1000$ and $K=200$ to problem \eqref{eqn: lasso} and display the comparison results in Figure \ref{fig: pgc vs existing methods}.
Analogously, one can see from this figure that the proposed PGC performs the best and outperforms the DLM and EXTRA.
Besides, the DLM, EXTRA and the proposed PGC all converge much faster than the DNG and DGD, which is consistent with the comparison results reported in \cite{shi14extra,Ling15decentralize}.
}

In Figure \ref{fig: spgc noise}, we present convergence curves of the proposed SPGC (Algorithm 1 with static graph), PG-EXTRA and D-SGD when the stochastic gradient information and the setting of Case 1 are used. The gradient noise power $\sigma^2$ is set to $0.1$ and $10$, respectively. 
Following Theorem \ref{thm:time:variant:inexact:rate}, the stepsize $\eta^r$ of the SPGC is set to $\eta_0 \sqrt{r}$, where $\eta_0\geq 0$. Note that when $\eta_0=0$, the SPGC reduces to the PGC in Algorithm 2. One can first observe from Figure \ref{fig: spgc noise}(a) that, due to the stochastic gradients, both the SPGC ($\sigma^2=0.1, \eta_0=2500$) and the PG-EXTRA ({$\sigma^2=0.1$}) suffer higher error floors than their counterparts with exact gradients in Figure \ref{fig: spgc no noise}. 
Moreover, it can also be seen from Figure \ref{fig: spgc noise}(b) that the PG-EXTRA ({$\sigma^2=0.1$})  achieves lower consensus error than the SPGC ($\sigma^2=0.1, \eta_0=2500$). 
However, as seen from Figure \ref{fig: spgc noise}(a), not only the SPGC ($\sigma^2=0.1, \eta_0=2500$) converges faster than the PG-EXTRA ({$\sigma^2=0.1$}), but also the achieved solution accuracy keeps decreasing with the iteration number. This is contrast to the PG-EXTRA ({$\sigma^2=0.1$}) whose accuracy is limited by an error floor. We can observe similar convergence results for the case with $\sigma^2=10$. Finally, it can be seen that the D-SGD converges much slower than the other methods. Note that the convergence curves of the D-SGD for $\sigma^2=0.1$ and $\sigma^2=10$ overlap, implying that the noisy gradients have less impact on the D-SGD.

In the last example, we examine the convergence behavior of the proposed DySPGC (Algorithm 1) over a time-varying network. Following Definition \ref{def:graph}, we assume that each link $(i,j)\in \cE$ has a probability $p_{ij}\in (0,1]$ being active (If $p_{ij}=1$ $\forall i,j$, then the DySPGC reduces to SPGC in static networks). 
The setting Case 1 is considered with gradient noise power $\sigma^2=0.1$, and the stepsize $\eta^r=2500\sqrt{r}$ of DySPGC is used. 
Figure \ref{fig: spgc time varying} displays the convergence curves of the DySPGC for various values of $p_{ij}$. As seen, the DySPGC exhibits considerable robustness against the time-varying networks.

\section{Conclusion}
In this paper we have proposed a dynamic stochastic proximal-gradient consensus (DySPGC) algorithm for solving a convex, possibly stochastic optimization problem over a randomly time-varying multi-agent network. We have analyzed the global convergence rate for DySPGC under various different scenarios, such as when the network is static/dynamic, or when the gradient is stochastic/deterministic. Our numerical results show that the proposed algorithms compare favorably with a few EXTRA based algorithms under various scenarios. 
Interestingly, Our algorithmic framework provides a unifying perspective for a few popular algorithms for distributed convex optimization. Such new perspective  allows significant generalization of these methods based on existing theories such as the primal-dual methods. As an example, leveraging upon the recent work \cite{Ouyang15}, we can develop an accelerated version of DySPGC, which is capable of further reducing certain constants in the convergence rate; see our technical report \cite{hong15consensus} for details.

There are a few interesting directions that we would like to pursue in the future. For example, can we generalize the algorithms and their analysis to problems with {\it nonconvex} objective functions?  Can we deal with a wider types of network dynamics such as the deterministic B-strongly connected networks? Is there a connection between distributed algorithms that we have studied in this work, with the optimization algorithms that minimizes a convex objective function consisting of a  finite sum of components (which often arise in machine learning related applications), such as the SAG algorithm \cite{Roux12astochastic} and the SVRG algorithm \cite{Johnson13}? {Some recent advancement in connecting distributed optimization methods with SAG/SVRG can be found in \cite{Mokhtari16DSA} and \cite{davood16nestt}}.

\vspace{-0.2cm}
\appendix

\subsection{Proof of Proposition \ref{prop:equiv2}}\label{app:implementation}
{The proof of this proposition is in fact straightforward. One only needs to start from the optimality condition of each step of Algorithm 2, and recognize that each  $\lambda^r$ has some symmetric structure, and that the iterates $\{z^r\}$ can be expressed using $\{x^r\}$.}

first let us split $\lambda^r$ by $\lambda^r=[\delta^r; \gamma^r]$ where $\delta^r, \gamma^r\in \mathbb{R}^{2EM\times 1}$.
The optimality conditions of \eqref{eq:x:update:static} -- \eqref{eq:lambda:update:static} are 
\begin{subequations}
\begin{align}
&G(x^{r}) +\zeta^{r+1}+A^T (\lambda^r + \Gamma(A x^{r+1}+ B z^r)) \nonumber\\
&\quad \quad \quad \quad\quad \quad\quad \quad\quad \quad+ \Omega (x^{r+1}-x^r) = 0\label{eq:opt:x}\\
& B^T (\lambda^r +  \Gamma(A x^{r+1}+ B z^{r+1})) = 0 \label{eq:opt:z}
\\
& \lambda^{r+1}-\lambda^r -\Gamma  \left( A x^{r+1}+ B z^{r+1}\right) = 0 \label{eq:opt:lambda}
\end{align}
\end{subequations}
where $\zeta^{r+1}\in\mathbb{R}^{NM}$ is a subgradient vector satisfying $\zeta^{r+1}_i\in h_i(x^{r+1}_i), \; \forall~i$.
First we show that
\begin{align}\label{eq:delta}
\delta^{r+1} = - \gamma^{r+1}, \quad\forall~r\ge 0.
\end{align}
At iteration $0$ this is true due to the initialization $B^T\lambda^0 = 0$. At iteration $r\ge 0$, by \eqref{eq:opt:z} and \eqref{eq:opt:lambda} we have
$B^T\lambda^{r+1}=0.$ This immediately implies
$
\delta^{r+1}_{ij}=-\gamma^{r+1}_{ij}, \; \forall~e_{ij}\in \cA.
$
Using this identity, we can define a new variable $\alpha$ as
$$\alpha^r= A^T\lambda^r = M_{-}\delta^r.$$
Applying \eqref{eq:delta} to \eqref{eq:opt:z} and the definition \eqref{eq:A}, we have
\begin{align}\label{eq:z:x:2}
z^{r+1} =\frac{1}{2} \left(A_1+A_2\right)x^{r+1} = \frac{1}{2} M^T_{+} x^{r+1},\; \forall~r\ge 0.
\end{align}
Using the above identity, we have
\begin{align*}
A x^{r+1} + B z^{r+1}&  
 = 
\frac{1}{2} \left[\begin{array}{l}
A_1-A_2 \\
A_2-A_1
\end{array}\right]x^{r+1}.
\end{align*}
This fact combined with the update rule for $\lambda$ implies that 
\begin{align}
\delta^{r+1}&= \delta^r+\frac{1}{2}(\Xi\otimes I_M) M^T_{-}x^{r+1}\nonumber\\
\alpha^{r+1} &
= \alpha^{r}+\frac{1}{2} M_{-}(\Xi\otimes I_M) M^T_{-} x^{r+1}.\label{eq:alaph:relation2}
\end{align}
By \eqref{eq:z:x:2}, we have 
\begin{align}
 A^T \Gamma B z^{r+1}=\frac{1}{2} M_{+}(\Xi\otimes I_M) M^T_{+}x^{r+1}, \quad \forall~r\ge 1 \nonumber
\end{align}
Utilizing the initial conditions $B^T \lambda^0=0$, $z^0=\frac{1}{2}M^T_{+}x^0$, and the fact that $A^T \lambda^{r+1}=\alpha^{r+1}$, the $x$-step optimality condition \eqref{eq:opt:x} can be written as
\begin{align}
&G(x^r)+ \zeta^{r+1} + \alpha^{r+1} + \frac{1}{2}M_{+}(\Xi\otimes I_M) M^T_{+} (x^{r+1}-x^r)\nonumber\\
&\quad \quad+ \Omega  (x^{r+1}-x^r) = 0\label{eq:opt:x:2}.
\end{align}
Plugging in \eqref{eq:M:Xi:product:b} and utilizing \eqref{eq:alaph:relation2}, \eqref{eq:opt:x:2} becomes
\begin{align}
&2\bigg(\sum_{j\in\cN_i}\frac{\rho_{ij}+\rho_{ji}}{2}+\frac{\omega_i}{2}\bigg) x^{r+1}_i  + \zeta^{r+1}\nonumber
\end{align}
\begin{align}
&\quad = -\nabla g_i(x_i^r) + \bigg(\sum_{j\in\cN_i}\frac{\rho_{ji}+\rho_{ij}}{2}+\frac{\omega_i}{2}\bigg)  x^r_i \nonumber
\\
&\quad +  \sum_{j\in \cN_i} \frac{\rho_{ji}+\rho_{ij}}{2} x^r_j +\frac{\omega_i}{2} x^r_i-\alpha^r_i, ~ \forall~i\label{eq:x:2}.
\end{align}
Moreover, \eqref{eq:alaph:relation2} can be expressed as, $\forall~i \in \cV$,
 \begin{align}
\alpha^{r+1}_i = \alpha^r_i +\sum_{j\in\cN_i}\frac{\rho_{ji}+\rho_{ij}}{2} x^{r+1}_i -\sum_{j\in\cN_i}\frac{\rho_{ji}+\rho_{ij}}{2} x^{r+1}_j.  \label{eq:alpha:2}
\end{align}
Next we remove the sequence $\{\alpha^{r}\}$ from the $x$ iterations. This is the key step towards obtaining a single-variable characterization.  To this end, we subtract \eqref{eq:x:2} by the same update for iteration $r$. By the definition of $\beta_i$ in \eqref{eq:beta}
we can derive the following update rule for node $i$ at iteration $r+1$:    
{\small
\begin{align}
&x^{r+1}_i-x^{r}_i  + \frac{1}{\beta_i}(\zeta_i^{r+1}-\zeta^r_i)-\frac{1}{\beta_i}\left(-\nabla g_i(x_i^r)+ \nabla g_i(x_i^{r-1})\right)\nonumber\\
&= \frac{1}{\beta_i} \sum_{j\in \cN_i} \hrho_{ij}(x^r_j-x^{r-1}_j)+(\frac{1}{2}+\frac{\omega_i}{2\beta_i})(x^r_i-x^{r-1}_{i}) \notag\\
&\quad +\frac{1}{\beta_i}(-\alpha^r_i+\alpha^{r-1}_i)\nonumber\\
&=   \frac{1}{\beta_i}  \sum_{j\in \cN_i} \hrho_{ij}(x^r_j-x^{r-1}_j)+(\frac{1}{2}+\frac{\omega_i}{2\beta_i})(x^r_i-x^{r-1}_{i})\nonumber\\
&\quad  -\frac{1}{\beta_{i}} \bigg(\sum_{j\in\cN_i}\hrho_{ij}x^{r}_i  -\sum_{j\in\cN_i}\hrho_{ij}x^{r}_j\bigg)~ (\text{by}~ \eqref{eq:alaph:relation2}, \eqref{eq:alpha:2})\nonumber\\
&= -  \frac{1}{2} x^{r-1}_i +\frac{\omega_i}{2\beta_i}(x^r_i-x^{r-1}_{i}) \notag \\
&\quad+ \frac{1}{2\sum_{j\in\cN_i}\hrho_{ij}+\omega_i} \sum_{j\in \cN_i} \hrho_{ij}(x^r_j-x^{r-1}_j)\nonumber\\ 
&\quad +\frac{1}{2\sum_{j\in\cN_i}\hrho_{ij}+\omega_i}\bigg(\sum_{j\in\cN_i}\hrho_{ij} x^{r}_j +\frac{\omega_i}{2}x^r_i\bigg)\nonumber
\\
&=  \frac{2}{\beta_i} \bigg(\sum_{j\in \cN_i} \hrho_{ij}x^r_j+\frac{\omega_i}{2}x^r_i\bigg) \notag \\
 &~~~~~~~-  \frac{1}{2} \bigg(x^{r-1}_i +  \frac{2}{\beta_i} \bigg(\sum_{j\in \cN_i}\hrho_{ij} x^{r-1}_j +\frac{\omega_i}{2}x^{r-1}_i\bigg)\bigg) \nonumber.
\end{align}}
Note that by the definition of $W$ in \eqref{eq:W}, we have {\small
\begin{align}\notag 
\frac{2}{\beta_i}\bigg(\sum_{j\in\cN_i}\hrho_{ij} x^{r}_j +\frac{\omega_i}{2} x^r_i \bigg)= \left[\left(W\otimes I_M\right) x^r\right]_{i}
\end{align}}
\!\!which is simply a weighted average of $x^r$ over all the neighbors of node $i$ (including itself).

Writing in vector form and utilizing the definition of $W$ in \eqref{eq:W}, we have
\begin{align*}
&x^{r+1}-x^r + \Upsilon^{-1}(\zeta^{r+1}-\zeta^r) - \Upsilon^{-1}\left(-G(x^r)+G(x^{r-1})\right)\nonumber\\
&\quad= (W\otimes I_M) x^r -\frac{1}{2}\left(I_{MN}+ W\otimes I_M\right) x^{r-1}, \quad\forall~r\ge 1.
\end{align*}
This proves the claim.

\vspace{-0.3cm}
\section{Preliminary Results}\label{app:prelim}
\subsection{Preliminary Results and Proof Outline}
{In this section we summarize a few preliminary results and identities that will be used later for proof of convergence of both Algorithm 1 and 2. We also provide a brief outline of the convergence proof. }

First we discuss the optimality condition for problem \eqref{eq:consensus:admm:matrix}.
Suppose that Assumption 1 holds. Let $y^*\in X^*$ denote an optimal solution of \eqref{eq:consensus}. Let $z^*_{ij}=y^*$, $(i,j)\in \cA$  and $x^*_i=y^*$, for all $i$. Due to equivalence of problems \eqref{eq:consensus} and \eqref{eq:consensus:admm}, $(z^*, x^*)$ is an optimal solution of \eqref{eq:consensus:admm:matrix}. From the assumed Slater condition we know that problem \eqref{eq:consensus:admm:matrix} has a saddle point $(x^*,z^*, \lambda^*)$ satisfying the following condition $\forall~z, \lambda,$ and $\forall~x \in\mbox{dom}(h)$
\begin{align}
L_0(x^*,z^*; \lambda)\le L_0(x^*,z^*; \lambda^*)\le & L_0(x,z; \lambda^*). \label{eq:saddle:1}
\end{align}
where $L_0(\cdot)$ is given by \eqref{eq:augmented} with $\Gamma\equiv 0$. 

Let us define the following vectors
\begin{align}
w:= [x; z; \lambda], \quad { F(w) := [A^T\lambda; B^T\lambda; -(Ax+Bz)].}\nonumber
\end{align}
The second inequality in \eqref{eq:saddle:1} is equivalent to the following
\begin{align}\label{eq:optimality2}
Q(w, w^*):=f(x)-f(x^*)+\langle w-w^*, F(w^*) \rangle\ge 0.
\end{align}
{The second inequality in \eqref{eq:saddle:1} also implies that $(x^*, z^*)$ is the optimizer for the following problem
	$$\min_{x\in \mbox{dom}(h), z} L_0(x, z; \lambda^*).$$
The first-order optimality condition of the above problem is given by the following (for some $\zeta^*\in\partial h(x^*)$)
\begin{align}\label{eq:optimality1}
\langle x-x^*, \nabla g(x^*)+\zeta^*\rangle + \langle w-w^*, F(w^*) \rangle \ge 0.
\end{align}
For notational simplicity, let us define the left hand side by $U(w, w^*)$. }

It is easy to observe that for all $x \in\mbox{dom}(h)$ and all $z, \lambda$, 
{\begin{align}\label{eq:F}
&\langle w-w^*, F(w^*) \rangle \nonumber\\
&=\langle x-x^*, A^T\lambda^* \rangle+ \langle y-y^*, B^T\lambda^* \rangle\nonumber\\ 
&=\langle x-x^*, A^T\lambda \rangle+ \langle y-y^*, B^T\lambda \rangle \nonumber\\
&\quad + \langle x-x^*, A^T(\lambda^* -\lambda) \rangle + \langle y-y^*, B^T(\lambda^* -\lambda) \rangle\nonumber\\ 
&=\langle x-x^*, A^T\lambda \rangle+ \langle y-y^*, B^T\lambda \rangle - \langle \lambda-\lambda^*, A x+ By\rangle \nonumber\\
&=\langle w-w^*, F(w) \rangle.
\end{align}
where in the second to the last equality we have used the fact that $Ax^*+ Bz^* =0$. }
Using the above identity, \eqref{eq:optimality2}--\eqref{eq:optimality1} are equivalent to the following two inequalities, respectively
\begin{subequations}
\begin{align}
\hspace{-0.3cm}Q(w,w^*)&=f(x)-f(x^*)+\langle w-w^*, F(w) \rangle\ge 0,\label{eq:optimality1:alter2}\\
\hspace{-0.3cm}U(w,w^*)&=\langle x-x^*, \nabla g(x^*)+\zeta^*\rangle+\langle w-w^*, F(w) \rangle\ge 0.\label{eq:optimality1:alter1}
\end{align}
\end{subequations}
Let us characterize the optimality condition for the iterates. Define a block diagonal matrix
${H}^{r+1}(\eta)$:
\begin{align*}
{H}^{r+1}(\eta):=
{\rm blkdg}\left\{\Omega+\eta I_{MN},(B^{r+1})^T \Gamma B^{r+1}, \Gamma^{-1} \right\}.
\end{align*}
{Let ${H}(\eta)$ denote its time-invariant counterpart, that is, replacing $B^{r+1}$ in the above definition by $B$}. Also define
$$F^{r}(w^r) := [(A^r)^T\lambda^r; (B^r)^T\lambda^r; -(A^r x^r+B^r z^r)].\nonumber$$
Using the fact that $\lambda^{r+1}=\lambda^r+\Gamma(Ax^{r+1}+Bz^{r+1})$, the optimality conditions for the subproblems of Algorithm 1 are given by [for all $x \in\mbox{dom}(h)$ and all $z, \lambda$]{\small
\begin{subequations}
\begin{align}
&\big\langle \tG^{r+1}(x^{r},\xi^{r+1}) +\zeta^{r+1} + (A^{r+1})^T\left(\lambda^{r+1}+\Gamma B^{r+1} (z^{r}-z^{r+1})\right)\nonumber\\
&~~~~+\left(\Omega + \eta^{r+1}I_{MN}\right)(x^{r+1}-x^r), x-x^{r+1} \big\rangle \ge 0, \label{eq:opt:x:admm2}\\
&\left\langle (B^{r+1})^T  \lambda^{r+1}, z-z^{r+1}\right\rangle \ge 0, \label{eq:opt:z:admm2}
\\
&\left\langle  \Gamma^{-1}(\lambda^{r+1}-\lambda^r) -\left( A^{r+1} x^{r+1}+ B^{r+1}z^{r+1}\right), \lambda-\lambda^{r+1}\right\rangle  \ge 0 \label{eq:opt:lambda:admm2}.
\end{align}
\end{subequations}}
Adding these conditions we obtain $\forall~x \in\mbox{dom}(h), \; \forall~ z,\lambda$, {
\begin{align}
&\left\langle x-x^{r+1}, \tG^{r+1}(x^{r},\xi^{r+1})+\zeta^{r+1}\right\rangle  \nonumber\\
&\quad +\langle w-w^{r+1}, F^{r+1}(w^{r+1}) \rangle +\langle (B^{r+1})^T \Gamma A^{r+1} (x-x^{r+1})  \notag \\
&\quad +(B^{r+1})^T \Gamma B^{r+1}(z-z^{r+1}), z^r- z^{r+1}\rangle\nonumber\\
&\quad - (w-w^{r+1})^T H^{r+1}(\eta^{r+1}) ({w}^{r}-w^{r+1})\ge 0.\nonumber
\end{align}}
Note that $(B^{r+1})^T \lambda^{r} = (B^r)^T\lambda^r =0$ because $\lambda^r=[\delta^{r}; -\delta^r]$, and each $B^{r+1}$ and $B^{r}$ stacks two identical matrices. Using this fact and the optimality condition of the $z$-step \eqref{eq:opt:z:admm2}, the following is true for any optimal solution $(z^*,x^*)$  
\begin{align}
&\big\langle (B^{r+1})^T \Gamma A^{r+1} (x^*-x^{r+1}) \nonumber\\
&\quad \quad +(B^{r+1})^T \Gamma B^{r+1}(z^*-z^{r+1}), z^r- z^{r+1}\big\rangle\le 0.\label{eq:z:opt:condition}
\end{align}
Combining the above two inequalities and rearranging terms, we obtain, for any $\tw:=(x^*; y^*; \lambda)$
\begin{align}\label{eq:admm:optimality}
&\left\langle x^*-x^{r+1}, G(x^{r})+\zeta^{r+1}\right\rangle  +\langle \tw-w^{r+1}, F^{r+1}(w^{r+1}) \rangle\notag \\
&~~~+\langle \tau^{r+1}, x^*-x^{r+1}\rangle \nonumber\\
&\ge  (\tw-w^{r+1})^T H^{r+1}(\eta^{r+1}) (w^{r}-w^{r+1}),
\end{align}
where we have defined the {\it gradient error} as
\begin{align}\label{eq:tau}
\tau^{r+1} :=\tG^{r+1}(x^{r},\xi^{r+1})-\nabla G(x^{r}).
\end{align}

{Next let us briefly provide an outline of the proof of convergence. To show convergence (i.e., Theorem \ref{tm:time:invariant}), we need to construct a {\it potential function} that decreases at each iteration. In our proof, the following quantity is used as the potential function
	\begin{align}\label{eq:potential}
	(w^* - w^{r+1})^T H(0) (w^* -w^{r+1}).
	\end{align}
	To show that such a measure decreases at each iteration, we need to utilize the optimality condition \eqref{eq:admm:optimality} that we have just derived from the execution of the algorithm, as well as the global optimality conditions \eqref{eq:optimality2} and \eqref{eq:optimality1}.

}

\vspace{-0.5cm}
\subsection{Proof of Theorem \ref{tm:time:invariant}}
We only prove the first part of the theorem. The second part is the consequence of Theorem \ref{thm:time:variant:inexact:rate}.
As $\eta^{r}=0$ for all $r$, and $\cG^r=\cG$ for all $r$, we denote $H:=H(0)$. Applying the static version of \eqref{eq:admm:optimality} and let $w^*:=(x^*, z^*, \lambda^*)$, we have
\begin{align}
&\left\langle x^*-x^{r+1}, G(x^r)+\zeta^{r+1}\right\rangle  +\langle w^*-w^{r+1}, F(w^{r+1}) \rangle \nonumber\\
&\ge  (w^*-w^{r+1})^T H(w^{r}-w^{r+1}).\nonumber
\end{align}
By using the convexity of $h$, we obtain{
\begin{align}\label{eq:convexity:h}
\left\langle x^*-x^{r+1}, \zeta^{r+1}\right\rangle &=  \left\langle x^*-x^{r+1}, \xi^{r+1}-\zeta^*+ \zeta^{*}\right\rangle\nonumber\\
&\le \left\langle x^*-x^{r+1},  \zeta^*\right\rangle,
\end{align}}
where $\zeta^*\in \partial h(x^{*})$. Similarly, we have
{\small\begin{align}\label{eq:convexity:g}
 &\langle x^* -x^{r+1}, G(x^r)\rangle\nonumber\\
 & = \langle x^* -x^{r+1}, G(x^r)-G(x^*)\rangle + \langle x^* -x^{r+1}, G(x^*)\rangle\nonumber\\
 &= \langle x^* -x^{r}, G(x^r)-G(x^*)\rangle + \langle x^r -x^{r+1}, G(x^r)-G(x^*)\rangle  \nonumber\\
 &\quad + \langle x^* -x^{r+1}, G(x^*)\rangle\nonumber\\
 &\stackrel{\rm(i)}\le -\|G(x^*)-G(x^r)\|^2_{\widetilde{P}^{-1}}+ \|x^r -x^{r+1}\|_{\widetilde{P}/4}^2\nonumber\\
 &\quad + \|G(x^r)-G(x^*)\|^2_{\widetilde{P}^{-1}} + \langle x^* -x^{r+1}, G(x^*)\rangle\nonumber\\
 &\le \|x^r -x^{r+1}\|_{\widetilde{P}/4}^2+ \langle x^* -x^{r+1}, G(x^*)\rangle
\end{align}}
\!\!where in $\rm (i)$ we have used the Young's inequality: $\langle a, b\rangle\le {\|a\|^2}/{(2\epsilon)} + {\epsilon\|b\|^2}/{2}$ for any $\epsilon>0$, and a key property due to Nesterov \cite[Theorem 2.1.5]{Nesterov04}. Namely, if $g_i(x_i)$ is convex with Lipschitzian gradient (constant $P_i$), then $\forall~ x_i, y_i\in X$,
\begin{align}
&\frac{1}{P_i}\|\nabla g_i(x_i)-\nabla g_i(y_i)\|^2 \le \langle \nabla g_i(x_i)-\nabla g_i(y_i), x_i-y_i\rangle.\nonumber
\end{align}
Combining the optimality condition \eqref{eq:optimality1:alter1} and the above two inequalities, we obtain
\begin{align}
\frac{1}{4}\|x^r-x^{r+1}\|_{\widetilde{P}}^2 \ge   (w^*-w^{r+1})^T H(w^{r}-w^{r+1}),\nonumber
\end{align}
which is equivalent to{
\begin{align}
&\|w^{*}-w^{r+1}\|^2_{H}\nonumber\\
&\le \frac{1}{2}\|x^r-x^{r+1}\|_{\widetilde{P}}^2 + \|w^{*}-w^{r}\|^2_{H} - \|w^{r}-w^{r+1}\|^2_H \label{eq:descent}.
\end{align}}
For time-invariant graph, \eqref{eq:z:x:2} is true,
which implies
\begin{align}
&(z^{r+1}-z^{r})^TB^T \Gamma B (z^{r+1}-z^{r})\nonumber\\
&=\frac{1}{4}(x^{r+1}-x^{r})^T M_{+} B^T \Gamma B M^T_{+}(x^{r+1}-x^{r}). \label{eq:relation:xz}
\end{align}
Plugging this relation into \eqref{eq:descent} we obtain
\begin{align}
&\|w^{*}-w^{r+1}\|^2_H \label{eq:condense}.\\
&\le \|w^{*}-w^{r}\|^2_H-\|x^r-x^{r+1}\|_{\Omega}-\Gamma^{-1} \|\lambda^r-\lambda^{r+1}\|^2\nonumber\\
&\quad - \frac{1}{4}(x^{r+1}-x^{r})^T \left({M_{+} B^T \Gamma B M^T_{+}}-2\widetilde{P}\right)(x^{r+1}-x^{r})  \nonumber.
\end{align}

Therefore, as long as
$\Omega+ \frac{1}{4}M_{+}B^T \Gamma B M^T_{+}-\frac{1}{2}\widetilde{P}\succ 0$
or equivalently
$
2\Omega+ M_{+}(\Xi \otimes I_{M}) M^T_{+}-\widetilde{P}\succ 0,
$  we will have $x^{r+1}\to x^r$, $\lambda^{r+1}\to \lambda^r$. By a standard argument (cf. the derivation in \cite[(A2.22)-(A2.25)]{chang14distributed}), we can argue that every limit point of the sequence $x^r$ and $\lambda^r$ is a primal dual optimal solution of problem \eqref{eq:consensus:admm:matrix}.
Finally, by noticing the identity
$2 \Omega + M_{+}(\Xi \otimes I_{M}) M^T_{+} = \Upsilon (W\otimes I_M+I_{MN})$
 by using the definitions of $W, \Upsilon$ and $M_{+}$, the theorem is proved.

\subsection{Proof of Theorem \ref{thm:time:variant:inexact:rate}}\label{app:thm:time:variant:inexact:rate}
Note that the assumption of boundedness of $x$ implies the boundedness of iterates $\{z^r\}$. This is because from the identity \eqref{eq:z:x:2} we have $z^{r+1}_{ij}=\frac{1}{2}(x^{r+1}_i+x^{r+1}_j)$ for all $r$. 
Therefore
{\small $$\|z^{r}-z^*\|_{B^T \Gamma B}^2 =\sum_{ij: e_{ij}\in
		\cA}2\rho_{ij}\bigg\|\frac{x^r_{i}+x^r_j}{2}-\frac{x^*_{i}+x^*_j}{2}\bigg\|^2\le d^2_z,\;\forall~r.$$}
\!\!By the convexity of $h$ and $g$ and the Lipschitz continuity of $\nabla g$, we have{\small
\begin{align*}
&\langle x^* -x^{r+1}, \zeta^{r+1}\rangle \le h(x^*) -h(x^{r+1})\\
&\langle x^* -x^{r+1}, G(x^r)\rangle = \langle x^* -x^r, G(x^r)\rangle +\langle x^r-x^{r+1}, G(x^r)\rangle\\
&\le g(x^*)-g(x^r) + g(x^r) -g(x^{r+1}) +\frac{1}{2}\|x^{r+1}-x^r\|^2_{\widetilde{P}}.
\end{align*}}
For any given $\lambda$, plugging in $\tw:=(x^*, z^*, \lambda)$ into \eqref{eq:admm:optimality} and use the previous two inequalities, we have
{	\begin{align}
	&{-Q({w}^{r+1}, \tw)}=f(x^*)-f(x^{r+1})+\langle \tw-w^{r+1}, F(w^{r+1})\rangle\nonumber\\
	&\ge   (\tw-w^{r+1})^T H(\eta^{r+1}) (w^{r}-w^{r+1})-\langle\tau^{r+1}, x^*-x^{r+1}\rangle\nonumber\\
	&\quad-\frac{1}{2}\|x^r-x^{r+1}\|_{\widetilde{P}}^2 \nonumber\\
	& \stackrel{\rm (i)}\ge (\tw-w^{r+1})^T H(\eta^{r+1}) (w^{r}-w^{r+1})-\langle\tau^{r+1}, x^*-x^{r}\rangle  \nonumber\\
	&\quad -\frac{\|\tau^{r+1}\|^2}{2\eta^{r+1}} -\frac{\eta^{r+1}}{2}\|x^{r+1}-x^{r}\|^2-\frac{1}{2}\|x^r-x^{r+1}\|_{\widetilde{P}}^2 \nonumber
		\end{align}
		\begin{align}
	&\ge \frac{1}{2}\|z^{r+1}-z^*\|_{B^T\Gamma B}^2-\frac{1}{2}\|z^{r}-z^*\|_{B^T\Gamma B}^2 + \frac{1}{2}\|\lambda^{r+1}-\lambda\|_{\Gamma^{-1}}^2\nonumber
	\\
	&\quad -\frac{1}{2}\|\lambda^{r}-\lambda\|_{\Gamma^{-1}}^2 -\langle\tau^{r+1}, x^*-x^{r}\rangle -\frac{\|\tau^{r+1}\|^2}{2\eta^{r+1}}\nonumber
		\\
	&+\frac{1}{2}\bigg(\|x^{r+1}-x^*\|_{\Omega+\eta^{r+1}I_{MN}}^2-\|x^{r}-x^*\|_{\Omega+\eta^{r+1}I_{MN}}^2\bigg),\label{eq:temp4}
	\end{align}}
\!\!where in $\rm (i)$ we have again used the Young's inequality; in the last inequality we have used the assumption \eqref{eq:Omega:inexact} (cf. the derivation in \eqref{eq:condense}).
Evaluating the LHS based on the average of the iterates $\barw^{r+1}$, and using convexity, we have {\small
	\begin{align}
	&{-Q(\bar{w}^{r+1}, \tw)}=f(x^*)-f(\bar{x}^{r+1})+\langle \tw-\barw^{r+1}, F(\barw^{r+1})\rangle\label{eq:Qin}\\
	&\ge \frac{1}{r+1}\sum_{t=0}^{r}\bigg(f(x^*)-f(x^{t+1})+\langle \tw-w^{t+1}, F(w^{t+1})\rangle\bigg)\nonumber\\
	&\ge -\frac{1}{2(r+1)}\|z^{0}-z^*\|_{B^T\Gamma B}^2- \frac{1}{2(r+1)}\|\lambda^0-\lambda\|_{\Gamma^{-1}}^2\nonumber\\
	&\quad + \frac{1}{2(r+1)}\sum_{t=0}^{r}\eta^{t+1}\bigg(\|x^{t+1}-x^*\|^2-\|x^{t}-x^*\|^2\bigg)\nonumber\\
	&\quad -\frac{1}{2(r+1)}\|x^{0}-x^*\|^2_{\Omega}-\frac{1}{r+1}\sum_{t=0}^{r}\bigg(\langle\tau^{t+1}, x^*-x^{t}\rangle +\frac{\|\tau^{t+1}\|^2}{2\eta^{t+1}}\bigg)\nonumber.
	\end{align}}
\!\!In the above derivation, the first inequality has used a similar identity as in \eqref{eq:F}, i.e.,
$\langle w-\tw, F(\tw) \rangle = \langle w-\tw, F(w) \rangle$. 
The rest of the proof follows the similar argument in the second part of the proof in \cite[Theorem 2.2]{gao14} [Eq. (25) -- Eq. (30)]. We include it here for completeness.
First let us take the supremum of $Q(\bar{w}^{r+1},\tw)$ over the ball $\mathcal{B}_{\rho}$. We have{
\begin{align}
&\sup_{\lambda\in\mathcal{B}_{\rho}}\left[ f(\barx^{r+1})-f(x^*)+\langle\bar{w}^{r+1}-\tw, F(\bar{w}^{r+1})\rangle\right]\nonumber\\
&=\sup_{\lambda\in\mathcal{B}_{\rho}}\bigg[ f(\barx^{r+1})-f(x^*)+\langle \barx^{r+1}-x^*, A^T\lambda\rangle\nonumber\\
&\quad+ \langle \bar{z}^{r+1}-z^*, B^T\lambda\rangle+ \langle \bar{\lambda}^{r+1}-\lambda, A\barx^{r+1}+B\bar{z}^{r+1}\rangle\bigg]\nonumber\\
&=\sup_{\lambda\in\mathcal{B}_{\rho}}\bigg[ f(\barx^{r+1})-f(x^*)+\langle \lambda, A\barx^{r+1}+B\bar{z}^{r+1}\rangle\bigg]\nonumber\\
&= f(\barx^{r+1})-f(x^*)+\rho\|A\barx^{r+1}+B\bar{z}^{r+1}\|. \label{eq:supQ}
\end{align}}
\!\!Further, we have the following series of inequalities{
\begin{align}
&\sum_{t=0}^{r}\sqrt{t+1}\left(\|x^{t}-x^*\|^2-\|x^{t+1}-x^*\|^2\right)\nonumber\\
&\le \sum_{t=0}^{r}\|x^{t+1}-x^*\|^2(\sqrt{t+1}-\sqrt{t}) \le d^2_x\sqrt{r+1}.
\end{align}}
\!\!Taking the supreme of both sides of \eqref{eq:Qin}, we obtain{
\begin{align}
&f(\barx^{r+1})-f(x^*)+\rho\|A\barx^{r+1}+B\bar{z}^{r+1}\|\nonumber\\
&\le \frac{1}{2(r+1)}\left(d^2_z +d^2_{\lambda}(\rho)+\max_i \omega_i d^2_x\right)\nonumber\\
&\quad +\frac{1}{2\sqrt{r+1}}d^2_x-\frac{1}{r+1}\sum_{t=0}^{r}\bigg(\langle\tau^{t+1}, x^*-x^{t}\rangle +\frac{\|\tau^{t+1}\|^2}{2\eta^{t+1}}\bigg)\notag 
\end{align}}
\!\!Taking the expectation on both sides of the above inequality and utilize the assumption made in \eqref{eq:gradient} about the stochastic gradient, and the fact that
{$$\sum_{t=0}^{r}\frac{1}{\sqrt{t+1}}\le \int_{1}^{r+1}\frac{1}{\sqrt{t}}dt\le 2\sqrt{r+1},$$}
we arrive at the desired bound.
\vspace{-0.3cm}

{
	\subsection{Proof of Theorem \ref{thm:random:exact}}
	Our proof is motivated by \cite{Wei13,chang14}. 
	Suppose at iteration $r$ we have iterate $w^r=(x^r, z^r, \lambda^r)$ and we are about to execute Algorithm 1. Consider the {\it virtual} sequence $(\hx^{r+1}, \hz^{r+1}, \hlambda^{r+1})$ generated by Algorithm 1 (based on $w^r$) with all nodes and edges being active (i.e., with $\cA^{r+1}=\cA$ and $\cV^{r+1}=\cV$). Then from \eqref{eq:descent} in the proof of Theorem \ref{tm:time:invariant}, we must have
	\begin{align}
	&\|w^{*}-\hw^{r+1}\|^2_H\le \frac{1}{2}\|x^r-\hx^{r+1}\|_{\widetilde{P}}^2 + \|w^{*}-w^{r}\|^2_H - \|w^{r}-\hw^{r+1}\|^2_H \label{eq:descent:invariant2}.
	\end{align}
	Define the following auxiliary sequences{
		\begin{align*}
		D_{x}(x^{r+1}, x^*): &= (x^*-x^{r+1})^T \Psi^{-1/2}\Omega \Psi^{-1/2}(x^*-x^{r+1})\\
		D_{z}(z^{r+1}, z^*): &= (z^*-z^{r+1})^T \Phi^{-1/2}B^T\Gamma B \Phi^{-1/2}(z^*-z^{r+1})\\
		D_{\lambda}(\lambda^{r+1}, \lambda^*): &= (\lambda^*-\lambda^{r+1})^T \Phi^{-1/2}\Gamma^{-1}\Phi^{-1/2}(\lambda^*-\lambda^{r+1}),
		\end{align*}}
	\!\!where $\Psi$ and $\Phi$ are given in \eqref{eq:Psi}. Also define $\cF^{r}=\{x^{t}, z^{t}, \lambda^t, \cG^t_d, t=1,\cdots, r\}$ as the filtration up to iteration $r$. The following is easy to verify{
		\begin{align}
		&\mathbb{E}\left[D_{x}(x^{r+1}, x^*)|\cF^{r}\right]  = \mathbb{E}\left[\sum_{i=1}^{N}\frac{\omega_i}{\alpha_i}\|x^{r+1}_i-x^*\|^2\right]\nonumber\\
		& = \sum_{i=1}^{N}w_i\|\hx^{r+1}_i-x^*\|^2 +\sum_{i=1}^{N} (1-\alpha_i)\frac{\omega_i}{\alpha_i}\|x^{r}_i-x^*\|^2\nonumber\\
		& = D_{x}(x^r, x^*)+ \|x^*-\hx^{r+1}\|^2_{\Omega} - \|x^*-x^{r}\|^2_{\Omega} \label{eq:exp:x}.
		\end{align}}
	Using the same argument, we have{
		\begin{align}
		&\mathbb{E}\left[D_{z}(z^{r+1}, z^*)|\cF^{r}\right] = D_{z}(z^r, z^*)+ \|z^*-\hz^{r+1}\|^2_{B^T\Gamma B} - \|z^*-z^{r}\|^2_{B^T\Gamma B} \label{eq:exp:z}\\
		&\mathbb{E}\left[D_{\lambda}(\lambda^{r+1}, \lambda^*)|\cF^{r}\right] = D_{\lambda}(\lambda^r, \lambda^*)+ \|\lambda^*-\hlambda^{r+1}\|^2_{\Gamma^{-1}} - \|\lambda^*-\lambda^{r}\|^2_{\Gamma^{-1}} \label{eq:exp:lambda}.
		\end{align}}
	Summing up \eqref{eq:exp:x} -- \eqref{eq:exp:lambda} and utilizing  \eqref{eq:descent:invariant2}, we obtain {
		\begin{align}
		&\mathbb{E}\left[D_{x}(x^{r+1}, x^*)|\cF^{r}\right]  + \mathbb{E}\left[D_{z}(z^{r+1}, z^*)|\cF^{r}\right] + \mathbb{E}\left[D_{\lambda}(\lambda^{r+1}, \lambda^*)|\cF^{r}\right]\nonumber\\
		&\le D_{x}(x^{r}, x^*)+D_{z}(z^{r}, z^*)+D_{\lambda}(\lambda^{r}, \lambda^*)+ \frac{1}{2}\|x^r-\hx^{r+1}\|_{\widetilde{P}}^2 \nonumber\\
		&\quad -(w^r-\hw^{r+1})^T H (w^{r}-\hw^{r+1})\nonumber\\
		&\le D_{x}(x^{r}, x^*)+D_{z}(z^{r}, z^*)+D_{\lambda}(\lambda^{r}, \lambda^*)- (\hx^{r+1}-x^r)^T\bigg({\Omega-\widetilde{P}/2}\bigg)(\hx^{r+1}-x^r) -\|\lambda^r-\hlambda^{r+1}\|^2_{\Gamma^{-1}} \nonumber
		\end{align}}
	where in the last inequality we have removed the term $(\hat{z}^{r+1}-z^r)^T B^T\Gamma B(\hat{z}^{r+1}-z^r)\succeq 0$.
	Using the assumption \eqref{eq:random:stepsize}, we conclude that the sequence
	$D_{x}(x^{r}, x^*)+D_{z}(z^{r}, z^*)+D_{\lambda}(\lambda^{r}, \lambda^*)$ is a nonnegative almost supermartingale, which is convergent by the nonnegative almost supermartigale convergence theorem \cite[Theorem 1]{robbins71}:
	{\begin{align*}
		&D_{x}(x^{r}, x^*), \; D_{z}(z^{r}, z^*),\; D_{\lambda}(\lambda^{r}, \lambda^*) \; \mbox{are bounded and converges w.p.1.}\\
		&\sum_{r=1}^{\infty}\|\lambda^r-\hlambda^{r+1}\|^2\le \infty, \quad \sum_{r=1}^{\infty}\|x^r-\hx^{r+1}\|^2\le \infty.
		\end{align*}}
	Then again by a standard argument (cf. \cite{Wei13}) we conclude that  $(x^r, z^r, \lambda^r)$ as well as $(\hx^r, \hz^r, \hlambda^r)$ converge with probability one to a primal-dual solution of problem \eqref{eq:consensus:admm:matrix}.

	\subsection{Proof of Theorem \ref{thm:random:inexact:rate}}
	Again we prove by utilizing the full iterates $\hw:=(\hx^{r+1}, \hz^{r+1}, \hlambda^{r+1})$. Define $\tw:= (x^*, z^*, \lambda)$ for any fixed $\lambda$.  From the derivation leading to \eqref{eq:temp4} we known that under the condition $\Omega\succ \widetilde{P}$, the full sequence satisfies{
		\begin{align}
		&-Q(\hw^{r+1}, \tw)=f(x^*)-f(\hx^{r+1})+\langle \tw-\hw^{r+1}, F(\hw^{r+1})\rangle\nonumber\\
		&\ge \frac{1}{2}\|\hz^{r+1}-z^*\|_{B^T\Gamma B}^2-\frac{1}{2}\|z^{r}-z^*\|_{B^T\Gamma B}^2 + \frac{1}{2}\|\hlambda^{r+1}-\lambda\|_{\Gamma^{-1}}^2\nonumber\\
		&+ \frac{1}{2}\left(\|\hx^{r+1}-x^*\|_{\Omega+\eta^{r+1}I_{MN}}^2-\|x^{r}-x^*\|_{\Omega+\eta^{r+1}I_{MN}}^2\right)\nonumber\\
		&-\frac{1}{2}\|\lambda^{r}-\lambda\|_{\Gamma^{-1}}^2-\langle\tau^{r+1}, x^*-x^{r}\rangle -\frac{\|\tau^{r+1}\|^2}{2\eta^{r+1}}. \label{eq:temp3}
		\end{align}}
	Notice that the following is true{
		\begin{align}
		&\langle \tw-\hw^{r+1}, F(\hw^{r+1})\rangle= \langle \tw-\hw^{r+1}, F(\tw)\rangle = -\langle \lambda, A\hx^{r+1}+B\hz^{r+1}\rangle.\label{eq:expression:F}
		\end{align}}
	\!\!Then it is easy to verify that {
		\begin{align}
		&\mathbb{E}\bigg[ \langle \tw-w^{r+1}, F(w^{r+1})\rangle   |\cF^{r}\bigg]= -\langle \lambda, A\Psi \hx^{r+1} + B\Phi \hz^{r+1}\rangle -\langle\lambda, Ax^r+Bx^r \rangle +\langle \lambda, A\Psi x^{r} + B\Phi z^{r}\rangle\label{eq:expectation:F}.
		\end{align}}
	Similarly, we have{\small
		\begin{align}
		\hspace{-0.3cm}\mathbb{E}\bigg[ \sum_{i=1}^{N}f_i(x^{r+1}_i)\big|\cF^{r}\bigg]=\sum_{i=1}^{N}\alpha_ i f_i(\hx^{r+1}_i)+\sum_{i=1}^{N}(1-\alpha_i)f_i(x^{r}_i)\label{eq:expectation:obj}.
		\end{align}}
	\!\!Using \eqref{eq:expectation:F}  -- \eqref{eq:expectation:obj} and the definition of $J(\cdot)$ in \eqref{eq:J}, the conditional expectation of $J(\cdot)$ can be expressed as below
	\begin{align}
	&\mathbb{E}[J(x^{r+1}, z^{r+1}, \lambda)|\cF^{r}] \nonumber\\
	&= J(x^{r}, z^{r}, \lambda) + \sum_{i=1}^{N} f_i(\hx^{r+1}_i) +\langle \lambda, A \hx^{r+1}+B\hz^{r+1}\rangle - \bigg(\sum_{i=1}^{N} f_i(x^{r}_i)+\langle \lambda, A x^{r}+Bz^{r}\rangle \bigg)\label{eq:long0}\\
	&\stackrel{\eqref{eq:expression:F}}=J(x^{r}, z^{r}, \lambda)+\sum_{i=1}^{N} f_i(\hx^{r+1}_i)-\langle \tw-\hw^{r+1}, F(\hw^{r+1})\rangle - \bigg(\sum_{i=1}^{N} f_i(x^{r}_i)-\langle \tw-w^{r}, F(w^{r})\rangle  \bigg)\nonumber\\
	&={ J(x^{r}, z^{r}, \lambda) + Q(\hw^{r+1},\tw)- Q(w^{r},\tw)}. \nonumber
	\end{align}

	Let us define{
		\begin{align}
		\widetilde{D}_{x}(x^{r+1}, x^*): &= \eta^{r+1}(x^*-x^{r+1})^T \Psi^{-1}(x^*-x^{r+1}),\nonumber
		\end{align}}
	then its conditional expectation is given by{
		\begin{align*}
		&\mathbb{E}\bigg[\widetilde{D}_{x}(x^{r+1}, x^*)|\cF^{r}\bigg]  = \mathbb{E}\bigg[\sum_{i=1}^{N}\frac{\eta^{r+1}}{\alpha_i}\|x^{r+1}_i-x^*\|^2\bigg]\nonumber\\
		&= \widetilde{D}_{x}(x^r, x^*)+ \|x^*-\hx^{r+1}\|_{\eta^{r+1}I_{MN}} - \|x^*-x^{r}\|^2_{\eta^{r+1}I_{MN}}+{(\sqrt{r+1}-\sqrt{r})\|x^{r}-x^*\|^2_{\Psi^{-1}}}.
		\end{align*}}
	Plugging \eqref{eq:temp3} and \eqref{eq:exp:x} -- \eqref{eq:exp:lambda} into \eqref{eq:long0}, we obtain a bound on the conditional expectation of $J(\cdot)$, given below
	\begin{align}
	&\mathbb{E}[J(x^{r+1}, z^{r+1}, \lambda)|\cF^{r}] \nonumber\\
	&\le J(x^{r}, z^{r}, \lambda) - Q(w^r, \tilde{w})+\frac{1}{2}\big(D_z(z^{r}, z^*)-\mathbb{E}[D_z(z^{r+1},z^*)]\big)+{(\sqrt{r+1}-\sqrt{r})\|x^r-x^*\|^2_{\Psi^{-1}}} \label{eq:long3}\\
	&+\frac{1}{2}\big(D_\lambda(\lambda^{r}, \lambda)-\mathbb{E}[D_\lambda(\lambda^{r+1},\lambda)]\big)+\frac{1}{2}\big(D_x(x^{r}, x^*)-\mathbb{E}[D_x(x^{r+1},x^*)]\big) \nonumber\\
	&\quad+ \frac{1}{2}\big(\widetilde{D}_x(x^{r}, x^*)-\mathbb{E}[\widetilde{D}_x(x^{r+1},x^*)]\big) + \langle\tau^{r+1}, x^*-x^{r}\rangle +\frac{\|\tau^{r+1}\|^2}{2\eta^{r+1}}.\nonumber
	\end{align}
	Let us define $\bar{x}^{r+1} = \frac{1}{r+1}\sum_{t=0}^{r}x^t$  and $\bar{z}^{r+1}$ similarly. Taking expectation wrt $\cF^r$ and summing over $t$, \eqref{eq:long3} becomes
	\begin{align}
	&{\mathbb{E}\left[Q(\bar{w}^r, \tw)\right]}\le \frac{1}{r+1}\mathbb{E}[J(x^0, z^0, \lambda)]+\frac{1}{2(r+1)} \mathbb{E}[D_z(z^0, z^*)]+\frac{1}{2(r+1)} \mathbb{E}[D_\lambda(\lambda^0, \lambda)] \nonumber\\
	&\quad+\frac{1}{r+1}\sum_{t=0}^{r} (\sqrt{t+1}-\sqrt{t})\|x^t-x^*\|^2_{\Psi^{-1}}+\frac{1}{2(r+1)} \mathbb{E}[D_x(x^0, x^*)]+\frac{1}{2(r+1)} \mathbb{E}[\widetilde{D}_x(x^0,x^*)]\nonumber\\
	&\quad+\frac{1}{r+1}\sum_{t=0}^{r}\left(\langle\tau^{t+1}, x^*-x^{t}\rangle +\frac{\|\tau^{t+1}\|^2}{2\eta^{t+1}}\right).\label{eq:long4}
	\end{align}

	By taking the superior of both sides over $\mathcal{B}_{\rho}:=\{\lambda\mid \|\lambda\|\le \rho\}$, and plugging in the definition of $Q(\bar{w}^r, \tw)$, we obtain
	\begin{align}
	&\mathbb{E}\big[\sum_{i=1}^{N}f_i(\bar{x}^r_i)-f(x^*)+\rho\|A\bar{x}^r+B\bar{y}^r\|\big]\nonumber\\
	& \le \frac{1}{r+1}\sup_{\lambda\in\mathcal{B}_{\rho}}J(x^0, z^0, \lambda) +\frac{1}{2(r+1)} D_z(z^0, z^*)+\frac{1}{2(r+1)} \sup_{\lambda\in\mathcal{B}_{\rho}}D_\lambda(\lambda^0, \lambda)+\frac{1}{r+1} D_x(x^0, x^*)\nonumber\\
	\quad&+\frac{1}{2(r+1)}\mathbb{E}[\widetilde{D}_x(x^0,x^*)]+\frac{1}{r+1}\sum_{t=0}^{r}\left(\langle\tau^{t+1}, x^*-x^{t}\rangle +\frac{\|\tau^{t+1}\|^2}{2\eta^{t+1}}\right)\nonumber\\
	&\quad+\frac{1}{r+1}\sum_{t=0}^{r} (\sqrt{r+1}-\sqrt{r})\|x^r-x^*\|^2_{\Psi^{-1}}\label{eq:long5}
	\end{align}

	The rest of the proof follows the last part of the proof of Theorem \ref{thm:time:variant:inexact:rate}.
	
	\subsection{Proof of Theorem \ref{thm:time:variant:inexact:rate:acc}}\label{app:acc}
	We first provide a lemma that bounds the quantity $Q(\cdot,\cdot)$ defined in \eqref{eq:optimality2}.
	\begin{lemma}\label{lm:Q}
		{\it  Let $\tw:=(x^*, y^*, \lambda)$ for any given $\lambda$. We have the following estimate for $Q(w^{r+1,\ag}, \tw)$
			\begin{align}
			&Q(w^{r+1,\ag}, \tw) - (1-\nu^r)Q(w^{r,\ag}, \tw)\nonumber\\
			&\le \nu^r\left(\langle \nabla g(x^{r+1,\md})+\zeta^*, x^{r+1}-x^*\rangle +\frac{\nu^r}{2}\|x^{r+1}-x^r\|^2_{\tP}+\langle w^{r+1}-\tw, F(\tw)\rangle\right)
			\end{align}
			for some $\zeta^*\in\partial h(x^*)$. }
		\begin{proof}
			From the definition of $x^{r+1,\ag}$, $x^{r+1,\md}$ we have
			\begin{align}\label{eq:agmd}
			x^{r+1,\ag}-x^{r+1,\md} = \nu^{r}(x^{r+1}-x^r).
			\end{align}
			First it is easy to show that for any feasible $x$, we have (cf. \cite[Eq. (2.16)]{Ouyang15})
			\begin{align}
			&g(x^{r+1,\ag}) - (1-\nu^{r}) g(x^{r,\ag}) \nonumber\\
			&\le \nu^r g(x^*)+\nu^{r}\langle G(x^{r+1,\md}), x^{r+1}-x^*\rangle+\frac{(\nu^r)^2}{2}\|x^{r+1}-x^{r+1}\|^2_{\tP}. \label{eq:g:acc}
			\end{align}
			Using this result, we have the following series of inequalities
			\begin{align}
			&Q(w^{r+1,\ag}, \tw) - (1-\nu^r)Q(w^{r,\ag}, \tw)\nonumber\\
			& = (f(x^{r+1,\ag})-f(x^*)) -(1-\nu^r) (f(x^{r,\ag})-f(x^*)) \nonumber\\
			&\quad \quad +\langle Ax^{r+1,\ag}+B z^{r+1,\ag},\lambda \rangle -(1-\nu^{r})\langle Ax^{r,\ag}+B z^{r,\ag},\lambda \rangle\nonumber\\
			& \le  \nu^{r}(h(x^{r+1})-h(x^*))+\nu^{r}\langle  G(x^{r+1,\md}), x^{r+1}-x^*\rangle\nonumber\\
			&\quad\quad+\frac{(\nu^r)^2}{2}\|x^{r+1}-x^{r+1}\|^2_{\tP}
			+\nu^r\langle Ax^{r+1}+Bz^{r+1},\lambda \rangle \nonumber
			\end{align}
			where the inequality uses \eqref{eq:g:acc}, the convexity of $h(\cdot)$ and the update rule of $x^{r+1,\ag}$.
		\end{proof}
	\end{lemma}
	
	We then proceed to prove Theorem \ref{thm:time:variant:inexact:rate:acc}.
	Let us define $\varpi^r = \frac{2}{(r+1)r}$. From \eqref{eq:nu} one can check that the following two identities hold
	$$\varpi^r = (1-\nu^r)\varpi^{r-1}, \quad \nu^r/\varpi^r=r, \quad \; \forall~r\ge 2.$$
	
	Let us define
	\begin{align*}
	{H}(\theta, \eta):=
	{\rm blkdg}\left\{\theta\Omega+\eta I_{MN},(B)^T \Gamma B, \Gamma^{-1} \right\}
	\end{align*}
	
	Similarly as in \eqref{eq:admm:optimality}, we can derive (for some $\zeta^{r+1}\in\partial h(x^{r+1})$)
	\begin{align*}
	&\left\langle x-x^{r+1}, G(x^{r+1,\md})+\zeta^{r+1}\right\rangle  +\langle w-w^{r+1}, F(w^{r+1}) \rangle +\langle \tau^{r+1}, x-x^{r+1}\rangle \nonumber\\
	&\ge  (w-w^{r+1})^T H(\theta^r, \eta^{r+1}) (w^{r}-w^{r+1}),\quad \forall~x\in\mbox{dom}(h), \; \forall~z,\lambda.
	\end{align*}
	Utilizing the fact that $\langle \tw-w^{r+1}, F(w^{r+1}) \rangle = \langle \tw-w^{r+1}, F(\tw) \rangle$ and that $\langle x^*-x^{r+1},\zeta^{r+1}-\zeta^*\rangle\le 0$ for any $\zeta^{r+1}\in\partial h(x^{r+1})$, we plugging in $\tw := (x^*, y^*, \lambda)$ and obtain 
	\begin{align*}
	&\left\langle x^*-x^{r+1}, G(x^{r+1,\md})+\zeta^{r+1}\right\rangle  +\langle \tilde{w}-w^{r+1}, F(\tw) \rangle +\langle \tau^{r+1}, x^*-x^{r+1}\rangle \nonumber\\
	&\ge  (\tilde{w}-w^{r+1})^T H(\theta^r, \eta^{r+1}) (w^{r}-w^{r+1}).
	\end{align*}
	Applying Lemma \ref{lm:Q}, we can obtain
	\begin{align}
	&Q(w^{r+1,\ag}, \tw) - (1-\nu^r)Q(w^{r,\ag}, \tw)\nonumber\\
	&\le \nu^r (\tw-w^{r+1})^T H(\theta^r, \eta^{r+1}) (w^{r+1}-w^{r}) +\frac{(\nu^r)^2}{2}\|x^{r+1}-x^r\|^2_{\tP}+ \nu^r\langle \tau^{r+1}, x^*-x^{r+1}\rangle\label{eq:temp5}.
	\end{align}
	From the assumption $\Omega$ should satisfy \eqref{eq:Omega:inexact:2}.
	Using such assumed bound, the definition of $\nu^r$ and $\theta^r$, and the fact that $\nu^r<1$ and $M_{+}(I_{M}\otimes\Xi) M^T_{+}\succeq 0$, we have
	\begin{align}
	\frac{\theta^r}{2}\Omega+ \frac{1}{4}M_{+}(I_{M}\otimes\Xi) M^T_{+}\succ \frac{\nu^r}{2}\widetilde{P}.
	\end{align}
	Applying the same derivation as in \eqref{eq:temp4}, and divide both sides of \eqref{eq:temp5} by $\varpi^r$,
	we obtain {
		\begin{align}
		&\frac{1}{\varpi^r}Q(w^{r+1,\ag}, \tw) - \frac{1-\nu^r}{\varpi^r}Q(w^{r,\ag}, \tw)\nonumber\\
		&=\frac{1}{\varpi^r}Q(w^{r+1,\ag}, \tw) - \frac{1}{\varpi^{r-1}}Q(w^{r,\ag}, \tw)\nonumber\\
		&\le \frac{\nu^r}{\varpi^r} \bigg((\tw-w^{r+1})^T H(\theta^r,\eta^{r+1}) (w^{r+1}-w^{r}) +\frac{\nu^r}{2}\|x^{r+1}-x^r\|^2_{\tP}+ \langle \tau^{r+1}, x^*-x^{r+1}\rangle\bigg)\nonumber\\
		&\le \frac{\nu^r}{2\varpi^r}\left(-\|z^{r+1}-z^*\|_{B^T\Gamma B}^2+\|z^{r}-z^*\|_{B^T\Gamma B}^2\right) + \frac{\nu^r}{2\varpi^r}\left(-\|\lambda^{r+1}+\lambda\|_{\Gamma^{-1}}^2+\|\lambda^{r}-\lambda\|_{\Gamma^{-1}}^2\right)\nonumber\\
		&-\quad \frac{\nu^r}{2 \varpi^r}\left(\|x^{r+1}-x^*\|_{\theta^r\Omega+\eta^{r+1}I_{MN}}^2+\|x^{r}-x^*\|_{\theta^r\Omega+\eta^{r+1}I_{MN}}^2\right)+\frac{\nu^r}{\varpi^r}\langle\tau^{r+1}, x^*-x^{r}\rangle +\frac{\nu^r}{\varpi^r}\frac{\|\tau^{r+1}\|^2}{2\eta^{r+1}}.
		\end{align}}
	Let us then analyze the successive sum of the RHS of the above inequality. Note that the sequences $\{\frac{\nu^r}{2\varpi^r}, \frac{\nu^r\eta^{r+1}}{2\varpi^r}\}$ are both increasing sequences, and the sequence $\frac{\nu^r\theta^r}{2\varpi^r}$ is non-increasing. Thus from \cite[Lemam 2.4]{Ouyang15} we have
	\begin{align}
	&\sum_{t=1}^{r}\frac{\nu^t}{2\varpi^t}\left(-\|z^{t+1}-z^*\|_{B^T\Gamma B}^2+\|z^{t}-z^*\|_{B^T\Gamma B}^2\right)\le \frac{\nu^r}{2\varpi^r} d^2_{z},\nonumber\\
	&\sum_{t=1}^{r}\frac{\nu^t}{2\varpi^t}\left(-\|\lambda^{t+1}+\lambda\|_{\Gamma^{-1}}^2+\|\lambda^{t}-\lambda\|_{\Gamma^{-1}}^2\right)\le \frac{\nu^r}{2\varpi^r} \sup_{\lambda^t}\|\lambda-\lambda^t\|_{\Gamma^{-1}}^2,\nonumber\\
	&\sum_{t=1}^{r}\frac{\nu^t\eta^{t+1}}{2\varpi^t}\left(-\|x^{t+1}-x^*\|^2+\|x^{t}-x^*\|^2\right)\le \frac{\nu^r\eta^r}{2\varpi^r} d^2_x\nonumber,\\
	&\sum_{t=1}^{r}\frac{\nu^t\theta^t}{2\varpi^t}\left(-\|x^{t+1}-x^*\|_{\Omega}^2+\|x^{t}-x^*\|_{\Omega}^2\right)\le \frac{\nu^1\theta^1}{2\varpi^1} \max_{i}\omega_i\|x^1-x^*\|^2 = \max_{i}\omega_i d^2_x \nonumber.
	\end{align}
	Combining these results, and noticing  $(1-\nu^1)/\varpi^1=0$, we obtain
	\begin{align}
	Q(w^{r+1,\ag}, \tw) &\le \frac{\nu^r}{2}\left(d^2_z+\sup_{\lambda^t}\|\lambda-\lambda^t\|_{\Gamma^{-1}}^2+\eta^r d^2_x\right)+\varpi^r \max_i d^2_x \omega_i\nonumber\\
	&\quad \quad + {\varpi^r}\sum_{t=1}^{r}\left(\frac{\nu^t}{\varpi^t}\frac{\|\tau^{t+1}\|^2}{2\eta^{t+1}}
	+\frac{\nu^t}{\varpi^t}\langle\tau^{t+1}, x^*-x^{t}\rangle\right).
	\end{align}
	Notice that from the derivation in \eqref{eq:supQ} we have
	\begin{align}
	\sup_{\lambda\in\mathcal{B}_{\rho}}Q(w^{r+1,\ag}, \tw) = f(x^{r+1,\ag})-f(x^*) +\rho\|Ax^{r+1,\ag}+Bz^{r+1,\ag}\|, \; \forall~\rho\ge 0.
	\end{align}
	Therefore we obtain 
	\begin{align}
	&f(x^{r+1,\ag})-f(x^*) +\rho\|Ax^{r+1,\ag}+Bz^{r+1,\ag}\|\nonumber\\
	&\le \frac{\nu^r}{2}\left(d^2_z+d^2_{\lambda}(\rho)+\eta^r d^2_x\right)+\varpi^r\max_i \omega_id^2_x+ {\varpi^r}\sum_{t=1}^{r}\left(\frac{\nu^t}{\varpi^t}\frac{\|\tau^{t+1}\|^2}{2\eta^{t+1}}+\frac{\nu^t}{\varpi^t}\langle
	\tau^{t+1}, x^*-x^{t}\rangle\right).
	\end{align}
	
	Taking expectation on both sides of the above inequality and utilizing
	\begin{align}
	&{\varpi^r}\sum_{t=1}^{r}\mathbb{E}\left[\frac{\nu^t}{\varpi^t}\frac{\|\tau^{t+1}\|^2}{2\eta^{t+1}}\right]=
	\sigma^2\varpi^r\sum_{t=1}^{r}\frac{\nu^t}{2\varpi^t \eta^{t+1}} = \frac{\sigma^2}{r(r+1)}\sum_{t=1}^{r}\frac{t}{\sqrt{t+1}}\le {\frac{2\sigma^2}{3}\frac{1}{\sqrt{r+1}}},\nonumber\\
	&\mathbb{E}\left[\frac{\nu^r}{\varpi^r}\langle\tau^{r+1}, x^*-x^{r}\rangle \right] =0, \;\forall~r.
	\end{align}
	we can obtain the desired bound.}

\bibliographystyle{IEEEbib}
\bibliography{ref,biblio,distributed_opt}

\end{document}